\magnification = 1100
\documentstyle{amsppt} 
\overfullrule0pt
\TagsOnLeft  
\def \n {\noindent}

\def \m {\medskip}
\def \b {\bigskip}
\def \a {\alpha}

\def \e {\epsilon} 
\def \End {{\operatorname {End}}}
  
\def\fsl{\frak s\frak l_2} 
\def\fsn{\frak s\frak l_n}
\def\fgl{\frak g\frak l_n}

\def \half {\frac{1}{2}}

\def \id {{\operatorname {Id}}}

\def \Hom {{\operatorname {Hom}}}
\def \l {\lambda}

\def \la {\langle}

\def \ot {\otimes}

\def \ra {\rangle}

\def \w {\omega}

\def \K {\Bbb K} 
\def \Z {\Bbb Z}
\font\bigf=cmbx10 scaled \magstep1

\document
\baselineskip=12pt 
\topmatter
\title {\bigf Representations of Two-Parameter \\
 Quantum Groups and Schur-Weyl Duality} \endtitle
\author Georgia Benkart \footnote {The authors gratefully acknowledge 
support from  National Science
Foundation Grant \#{}DMS--9970119, National Security Agency Grant
\#{}MDA904-01-1-0067, and the hospitality of
the Mathematical Sciences Research Institute, Berkeley. \ \ \ \ \ \ \
 \ \ \ \ \ \ \ \ \ \ } \\ 
Sarah Witherspoon$^1$ \endauthor
\leftheadtext{ GEORGIA BENKART, \; SARAH WITHERSPOON} 
\rightheadtext {TWO-PARAMETER QUANTUM GROUPS} 
\m
\abstract \quad We determine the finite-dimensional simple modules for 
two-parameter quantum groups corresponding to the general linear and 
special linear Lie algebras ${\frak {gl}}_n$ and ${\frak {sl}}_n$, and 
give a complete reducibility result.  These quantum groups have a 
natural $n$-dimensional module $V$.  We prove an analogue of 
Schur-Weyl duality in this setting: \ the centralizer algebra of the 
quantum group action on the $k$-fold tensor power of $V$ is a quotient 
of a Hecke algebra for all $n$ and is isomorphic to the Hecke algebra in case 
$n\geq k$.  \endabstract \date August 3, 2001 \enddate \subjclass\nofrills 
\n 2000 {\it Mathematics Subject Classification.\/} Primary 17B37, 
16W30, 16W35, 81R50 \endsubjclass

\endtopmatter 

 \head Introduction \endhead \m In this work we 
study the representations of two two-parameter quantum groups 
$\widetilde{U}=U_{r,s}({\frak {gl}}_n)$ and $U=U_{r,s}({\frak 
{sl}}_n)$.  Our Hopf algebra $\widetilde{U}$ is isomorphic as an 
algebra to Takeuchi's $U_{r,s^{-1}}$ (see [T]), but as a Hopf algebra, 
it has the opposite coproduct.  As an algebra, $\widetilde U$ has generators $e_j, \ f_j, (1 
\leq j < n)$, and $\ a_i^{\pm 1}, \ b_i^{\pm 1}$ \ \ ($1 \leq i \leq 
n$), and defining relations given in (R1)-(R7) below. 
The elements $\ a_i^{\pm 1}, \ b_i^{\pm 1}$ generate a commutative 
subalgebra $\widetilde U^0$, and the elements $e_j, \ f_j$, \ 
$\w_j^{\pm 1}, (\w_j')^{\pm 1}$ ($1 \leq j < n)$, where $\w_j = 
a_jb_{j+1}$ and $\w_j' = a_{j+1}b_j$, generate the subalgebra $U = 
U_{r,s}(\fsn)$.
\m
 The structure of these quantum groups 
was investigated in [BW], where we realized both $\widetilde{U}$ and 
$U$ as Drinfel'd doubles of certain Hopf subalgebras and constructed 
an $R$-matrix for $\widetilde{U}$ and $U$.  In particular, for any two 
$\widetilde{U}$-modules in category $\Cal O$ (defined in Section 3), there is an 
isomorphism $R_{M',M}: M'\otimes M\rightarrow M\otimes M'$.  The 
construction of $R_{M',M}$ is summarized in Section 4 of this note.  
In Sections 2 and 3, we classify the finite-dimensional simple 
$\widetilde{U}$-modules when $rs^{-1}$ is not a root of unity and 
prove that all finite-dimensional $\widetilde{U}$-modules on which 
$\widetilde{U}^0$ acts semisimply are completely reducible.  These 
results hold equally well for $U$.  The hypothesis on 
$\widetilde{U}^0$ is necessary: we provide examples of 
finite-dimensional modules that are not completely reducible.  Our 
complete reducibility proof uses a quantum Casimir operator defined in 
[BW] and parallels the argument in [L].  \m There is a natural 
$n$-dimensional module $V$ for $\widetilde U$ (resp., $U$) defined in 
Section 1.  On tensor powers $V^{\ot k}$ of $V$, the transformations 
$R_i = \id^{\ot(i-1)} \ot R_{V,V} \ot \id^{\ot (k-i-1)}$ ($1 \leq i < 
k$) commute with the action of $\widetilde U$, and so they generate a 
subalgebra of $\End_{\widetilde U}(V^{\ot k})$.  This yields a map 
from a two-parameter Hecke algebra $H_k(r,s)$ to 
$\End_{\widetilde{U}}(V^{\ot k})$.  In the final section we prove a 
two-parameter analogue of Schur-Weyl duality: \ The transformations 
$R_i$ generate the {\it full} centralizer algebra 
$\End_{\widetilde{U}}(V^{\otimes k})$, and in case $n\geq k$, this 
centralizer algebra is isomorphic to $H_k(r,s)$.  The proof is 
elementary relying only on basic facts about the representations and 
explicit computations, and we believe it is new in the one-parameter 
case as well (compare [Ji]).  It is similar to a proof of classical 
Schur-Weyl duality due to the first author, which can be found in [H].  
An interesting consequence of the argument is the result that $V^{\ot 
k}$ is a cyclic $\widetilde U$-module for $n \geq k$.  \m Throughout 
we will work over an algebraically closed field $\K$.  \m \b \head \S 
1.  Preliminaries \endhead \m First we recall the definitions of the 
two-parameter quantum groups from [BW], and some basics about their 
representations.  Let $\e_1, \dots, \e_n$ denote an orthonormal basis 
of a Euclidean space $E$ with an inner product $\langle \ , \ 
\rangle$.  Let $\Pi = \{\a_j = \e_{j}- \e_{j+1} \mid j = 1, \dots, 
n-1\}$ and $\Phi = \{\e_i -\e_j \mid 1 \leq i \neq j \leq n\}$.  Then 
$\Phi$ is a finite root system of type A$_{n-1}$ with $\Pi$ a base of 
simple roots.  \m Fix nonzero elements $r,s$ in $\K$ with $r \neq s$.  
\m Let $\widetilde U = U_{r,s}(\fgl)$ be the unital associative 
algebra over $\K$ generated by elements $e_j, \ f_j, (1 \leq j < n)$, 
and $\ a_i^{\pm 1}, \ b_i^{\pm 1}$ \ \ ($1 \leq i \leq n$), which 
satisfy the following relations.  \m \roster \item"{(R1)}" The 
$a_i^{\pm 1}, \ b_j^{\pm 1}$ all commute with one another and \ 
$a_ia_i^{-1}= b_j b_j^{-1}=1,$ \m \item"{(R2)}" $ a_ie_j = 
r^{\la\epsilon_i,\a_j\ra}e_j a_i$ \ \ and \ \ $a_if_j 
=r^{-\la\epsilon_i,\a_j\ra} f_ja_i,$ \m \item"{(R3)}" $b_ie_j = 
s^{\la\epsilon_i,\a_j\ra}e_j b_i$ \ \ and \ \ $b_if_j =s^{-\la 
\epsilon_i,\a_j\ra} f_jb_i,$ \m \item"{(R4)}" 
$\displaystyle{[e_i,f_j]=\frac{\delta_{i,j}}{r-s}(a_ib_{i+1}-a_{i+1}b_i),}$ 
\m \item"{(R5)}" $[e_i,e_j]=[f_i,f_j]=0 \ \ \text{ if }\ \ |i-j|>1, $ 
\m \item"{(R6)}" $e_i^2e_{i+1}-(r+s)e_ie_{i+1}e_i+rse_{i+1}e_i^2 = 0,$ 
\smallskip \hskip -.2 truein $e_i e^2_{i+1} -(r+s)e_{i+1}e_ie_{i+1} 
+rs e^2_{i+1}e_i = 0,$ \m \item"{(R7)}" 
$f_i^2f_{i+1}-(r^{-1}+s^{-1})f_if_{i+1}f_i +r^{-1}s^{-1}f_{i+1}f_i^2 = 
0,$ \smallskip \hskip -.2 truein $f_i f^2_{i+1} 
-(r^{-1}+s^{-1})f_{i+1}f_if_{i+1}+r^{-1}s^{-1} f^2_{i+1} f_i=0.$ 
\endroster \m We will be interested in the subalgebra $U = 
U_{r,s}(\fsn)$ of $\widetilde U = U_{r,s}(\fgl)$ generated by the 
elements $e_j,f_j$, $\w_j$, and $\w_j'$ \ ($1 \leq j < n)$, where

$$\omega_j = a_j b_{j+1}  \ \ \text{and} \ \ \omega_j' = a_{j+1}b_j. \tag 1.1$$
\n 
These elements satisfy (R5)-(R7) along with the following relations:
\m
\roster
\item"{(R1')}"  The $\w_i^{\pm 1}, \ \w_j^{\pm 1}$ all commute with one 
another and $\w_i \w_i^{-1}= \w_j'(\w_j')^{-1}=1,$
\m
\item"{(R2')}"  $ \w_i e_j = r^{\la \e_i,\a_j\ra }s^{\la \e_{i+1},\a_j\ra}e_j \w_i$ \ \ and
\ \
$\w_if_j = r^{-\la\e_i,\a_j\ra}s^{-\la\e_{i+1},\a_j\ra} f_j\w_i,$ 
\m
\item"{(R3')}"
$\w_i'e_j = r^{\la\epsilon_{i+1},\a_j\ra}s^{\la\e_{i},\a_j\ra}e_j \w_i'$  \ \ and
\ \
$\w_i'f_j = r^{-\la\epsilon_{i+1},\a_j\ra}s^{-\la\e_{i},\a_j\ra}f_j \w_i'$,  
\m 
\item"{(R4')}"
$\displaystyle{[e_i,f_j]=\frac{\delta_{i,j}}{r-s}(\w_i-\w_i').}$
\endroster
\m 
When $r = q$ and $s = q^{-1}$, the algebra $U_{r,s}(\fgl)$ modulo
the ideal generated by the elements $b_i-a_i^{-1}$, $1 \leq i \leq n$,  is just the
quantum general linear group $U_q(\fgl)$,  and $U_{r,s}(\fsn)$ modulo
the ideal generated by the elements $\w_j'-\w_j^{-1}$, $1 \leq j < n$, is $U_q(\fsn)$.
\m
 The algebras $\widetilde U$ and $U$ are Hopf algebras, where the 
$a_i^{\pm 1}, b_i^{\pm 1}$ are group-like elements, and the remaining 
Hopf structure is given by

$$\aligned
& \Delta(e_i) =e_i\otimes 1 + \omega_i\otimes e_i, \quad  \ 
\Delta(f_i)=1\otimes f_i + f_i\otimes \omega_i',\\
& \varepsilon(e_i)=\varepsilon(f_i)=0, \ \ \
\, S(e_i)=-\omega_i^{-1}e_i, \ \ \ S(f_i)= -f_i(\omega_i')^{-1}.
\endaligned \tag 1.2$$

\m
Let $\Lambda= \Z \e_1 \oplus \cdots \oplus \Z \e_n$, the weight
lattice of $\fgl$, and $Q=\Z \Phi$ the root lattice.    
Corresponding to any  $\l \in \Lambda$
is an algebra homomorphism $\hat \l$ from the subalgebra
$\widetilde U^{0}$ of $\widetilde U$ generated by the elements $a_i^{\pm
1}, b_i^{\pm 1}$ ($1 \leq i \leq n$) to $\K$ given by 

$$ \hat \l(a_i) =  r^{\la\e_i,\lambda\ra} \qquad \text{and} \qquad 
\hat \l(b_i) = s^{\la\e_{i},\lambda\ra}. \tag 1.3$$ 

\n The restriction $\hat \l: U^{0} \rightarrow \K$ of $\hat{\lambda}$ to
the subalgebra $U^0$ of
$U$ generated by $\w_j^{\pm 1}, (\w_j')^{\pm 1}$ ($1 \leq j < n$)
satisfies

$$\hat \l(\w_j) = r^{\la\e_j,\l\ra}s^{\la\e_{j+1},\l\ra} \qquad \text{and} \qquad 
\hat \l(\w_j') = r^{\la\e_{j+1},\l\ra}s^{\la\e_j,\l\ra}. \tag 1.4$$ 
\m 
Let $M$ be a module for $\widetilde U = U_{r,s}(\fgl)$ of dimension $d < \infty$.
As $\K$ is algebraically closed, we have

$$M = \bigoplus_{\chi} M_\chi,$$

\n where each $\chi: \widetilde U^0 \rightarrow \K$ is an algebra homomorphism,
and  $M_\chi$ is the generalized eigenspace given by

$$M_\chi = \{m \in M \mid (a_i-\chi(a_i)\,1)^d m
= 0 =  (b_i-\chi(b_i)\,1)^d m, \ \ \text{ for all} \ i\}. \tag 1.5$$

\n  When $M_\chi \neq 0$ we say that $\chi$ is
a {\it weight} and $M_{\chi}$ is the corresponding {\it weight space}. 
(If $M$ decomposes into genuine eigenspaces relative
to $\widetilde U^0$ (resp. $U^0$), then we say that
$\widetilde U^0$ (resp. $U^0$) {\it acts semisimply on} $M$.)   
\m From relations (R2) and (R3) we deduce  that 

$$\aligned
& e_j M_\chi \subseteq M_{\chi \cdot \widehat {\a_j}} \\
& f_j M_\chi \subseteq M_{\chi \cdot (\widehat {-\a_j})}\endaligned \tag 1.6$$

\n where $\widehat {\a_j}$ is as in (1.3), and $\chi \cdot \psi$ is the 
homomorphism with values $(\chi \cdot \psi)(a_i) = \chi(a_i)\psi(a_i)$ 
\ and \ $(\chi \cdot \psi)(b_i) = \chi(b_i)\psi(b_i)$.  In fact, if 
$(a_i-\chi(a_i)\,1)^k m = 0$, then applying relation (R2) yields $(a_i 
- \chi(a_i) r^{\la\e_i,\a_j\ra}\,1)^k e_j m = 0$, and similarly for 
$b_i$ and for $f_j$.  Therefore, the sum of eigenspaces is a submodule 
of $M$, and if $M$ is simple this sum must be $M$ itself.  Thus, in 
(1.5), we may replace the power $d$ by 1 whenever $M$ is simple, and 
$\widetilde U^0$ must act semisimply in this case.  We also can see 
from (1.6) that for each simple $M$ there is a homomorphism $\chi$ so 
that all the weights of $M$ are of the form $\chi \cdot \hat \zeta$, 
where $\zeta \in Q$.
\m
It is shown in [BW, Prop.  3.5] that if $\hat{\zeta}=\hat{\eta}$, then 
$\zeta=\eta \ (\zeta,\eta\in Q)$ provided
$rs^{-1}$ is not a root of unity.  As a result, we have the 
following proposition.  \b \proclaim{Proposition 1.7} [BW, Cor.\ 3.14] 
Let $M$ be a finite-dimensional module for $U_{r,s}({\frak {sl}}_n)$ 
or for $U_{r,s}({\frak {gl}}_n)$.  If $rs^{-1}$ is not a root of 
unity, then the elements $e_i,f_i$ $(1\leq i < n)$ act nilpotently on 
$M$.  \endproclaim \b When $rs^{-1}$ is not a root of unity, a 
finite-dimensional simple module $M$ is a {\it highest weight} module 
by Proposition 1.7 and (1.6).  Thus there is some weight $\psi$ and a 
nonzero vector $v_0\in M_{\psi}$ such that $e_jv_0=0$ for all 
$j=1,\ldots, n-1$, and $M=\widetilde{U}.v_0$.  It follows from the 
defining relations that $\widetilde{U}$ has a triangular 
decomposition: $\widetilde{U}=U^-\widetilde{U}^0U^+$, where $U^+$ 
(resp., $U^-$) is the subalgebra generated by the elements $e_i$ 
(resp., $f_i$).  Applying this decomposition to $v_0$, we see that 
$M=\oplus_{\zeta\in Q^+} M_{\psi\cdot (\widehat{-\zeta})}$, where 
$Q^+=\sum_{i=1}^{n-1}\Z_{\geq 0}\alpha_i$.  \m When all the weights of 
a module $M$ are of the form $\hat{\lambda}$, where $\lambda\in 
\Lambda$, then for brevity we say that $M$ has weights in $\Lambda$.  
Rather than writing $M_{\hat{\lambda}}$ for the weight space, we 
simplify the notation by writing $M_{\lambda}$.  Note then (1.6) can 
be rewritten as $e_jM_{\lambda}\subseteq M_{\lambda +\alpha_j}$ and 
$f_jM_{\lambda}\subseteq M_{\lambda -\alpha_j}$.  Any simple 
$\widetilde{U}$-module having one weight in $\Lambda$ has all its 
weights in $\Lambda$.  \m Next we give an example of a simple 
$\widetilde{U}$-module with weights in $\Lambda$, which is the 
analogue of the natural representation for ${\frak {gl}}_n$.  \b 
\subhead The natural representation for $U_{r,s}(\fgl)$ and 
$U_{r,s}(\fsn)$ \endsubhead \m Consider an $n$-dimensional vector 
space $V$ over $\K$ with basis $\{v_j \mid 1 \leq j \leq n\}$.  We 
define an action of the generators of $\widetilde U = U_{r,s}(\fgl)$ 
by specifying their matrices relative to this basis:

$$\gathered
e_j = E_{j,j+1}, \qquad \  f_j = E_{j+1,j}, \qquad (1 \leq j < n) \\
a_i = rE_{i,i}+ \sum_{k \neq i} E_{k,k}, \qquad \qquad \quad \  (1 \leq i \leq n) \\
b_i = sE_{i,i} + \sum_{k \neq i} E_{k,k} \qquad \qquad \qquad  (1 \leq i \leq n).
\endgathered $$  

\n It follows that  $\w_j = a_jb_{j+1} =  rE_{j,j} + sE_{j+1,j+1} + \sum_{k \neq j,j+1}
E_{k,k}$ and $\w_j' = a_{j+1}b_j = sE_{j,j} + rE_{j+1,j+1} +  \sum_{k \neq j,j+1}
E_{k,k}$.
Now to verify that this extends to an action of
$\widetilde U$, (hence of $U = U_{r,s}(\fsn)$),  we need to check that the
relations hold.  We present an illustrative example and
leave the remainder to the reader: $$\aligned a_i e_j& = (rE_{i,i} + 
\sum_{k \neq i} E_{k,k})E_{j,j+1} \\
& = \cases
rE_{j,j+1}  & \quad \text{if $j=i$} \\
E_{j,j+1}  & \quad \text{if $j\neq i$}. 
\endcases \endaligned $$ 

\n This can be seen to equal $r^{\la\e_i,\a_j\ra}E_{j,j+1}(rE_{i,i} +
\sum_{k \neq i} E_{k,k})$,
which confirms that 
$a_i e_j = r^{\la\e_i,\a_j\ra} e_j a_i$
holds.   
\m
It follows from the fact that $a_i v_j = r^{\la\e_i,\e_j\ra} v_j$ and
$b_i v_j = s^{\la\e_i,\e_j\ra} v_j$ for all $i,j$ 
that $v_j$ corresponds to the weight $\e_j
= \e_1 -(\a_1 + \cdots + \a_{j-1})$.  
Thus, $V = \bigoplus_{j = 1}^n V_{{\e_j}}$ is the natural analogue of 
the $n$-dimensional representation of $\fgl$ and $\fsn$, and it is a 
simple module for both $\widetilde U$ and $U$.  When $r=q$ and 
$s=q^{-1}$, $b_i$ acts as $a_i^{-1}$ on $V$, and so $V$ is a module 
for the quotient $U_q({\fgl})$ of $U_{q,q^{-1}}({\fgl})$ by the ideal 
generated by $b_i-a_i^{-1} \ (1\leq i\leq n)$.  This is the natural 
module for the one-parameter quantum group $U_q({\fgl})$.  A similar 
statement is true for $U_q({\frak {sl}}_n)$.  \b \m \head \S 2.  
Classification of finite-dimensional simple modules \endhead \m 
Results will be stated for $\widetilde{U}$-modules, but everything 
holds as well for $U$-modules.  \m Let $\widetilde{U}^{\geq 0}$ denote 
the subalgebra of $\widetilde{U}$ generated by $a_i,b_i \ (1\leq i 
\leq n)$ and $e_i \ (1\leq i <n)$.  Let $\psi$ be any algebra 
homomorphism from $\widetilde{U}^0$ to $\K$, and $V^{\psi}$ be the 
one-dimensional $\widetilde{U}^{\geq 0}$-module on which $e_i$ acts as 
multiplication by 0 $(1\leq i < n)$, and $\widetilde{U}^0$ acts via 
$\psi$.  We define the {\it Verma module} $M(\psi)$ with highest 
weight $\psi$ to be the $\widetilde{U}$-module induced from 
$V^{\psi}$, that is
$$M(\psi)= \widetilde{U}\otimes _{\widetilde{U}^{\geq 0}}V^{\psi}.$$
Let $v_{\psi}=1\otimes v\in M(\psi)$, where $v$ is any nonzero vector
of $V^{\psi}$.  Then $e_i . v_{\psi}=0 \ (1\leq i <n)$ and $a . v_{\psi}
=\chi(a)v_{\psi}$ for any $a\in \widetilde{U}^0$ by construction.
\m
Notice that $\widetilde{U}^0$ acts semisimply on $M(\psi)$ by relations
(R2) and (R3).  If $N$
is a $\widetilde{U}$-submodule of $M(\psi)$, then $N$ is also a 
$\widetilde{U}^0$-submodule of the $\widetilde{U}^0$-module $M(\psi)$,
and so $\widetilde{U}^0$ acts semisimply on $N$ as well.  If $N$ is a
{\it proper} submodule, it must be that $N\subset \sum_{\mu\in Q^+\setminus 
\{0\}} M(\psi)_{\psi\cdot (\widehat{-\mu})}$ by (1.6), as 
$M(\psi)_{\psi}=\K v_{\psi}$ generates $M(\psi)$.  Therefore $M(\psi)$ 
has a unique maximal submodule, namely the sum of all proper 
submodules, and a unique
simple quotient, $L(\psi)$.    In fact, all 
finite-dimensional simple $\widetilde{U}$-modules are of this form, as 
the following theorem demonstrates.  \b \proclaim{Theorem 2.1} Let 
$\psi:\widetilde{U}^0\rightarrow \K$ be an algebra homomorphism.  Let 
$M$ be a $\widetilde{U}$-module, on which $\widetilde{U}^0$ acts 
semisimply and which contains an element $m\in M_{\psi}$ such that $e_i .  
m=0$ for all $i$ $(1\leq i <n)$.  Then there is a unique homomorphism 
of $\widetilde{U}$-modules $F:M(\psi)\rightarrow M$ with $F(v_{\psi})= 
m$.  In particular, if $rs^{-1}$ is not a root of unity and $M$ is a 
finite-dimensional simple $\widetilde{U}$-module, then $M\cong 
L(\psi)$ for some weight $\psi$.  \endproclaim \m \demo{Proof} By the 
hypothesis on $m$, ${\K} m$ is a one-dimensional $\widetilde{U}^{\geq 
0}$-submodule of $M$, considered as a $\widetilde{U} ^{\geq 0}$-module 
by restriction.  In fact, mapping $v_{\psi}$ to $m$ yields a 
$\widetilde{U}^{\geq 0}$-homomorphism from $V^{\psi}$ to $\K m$.  By 
the definition of $M(\psi)$, we have 
$\Hom_{\widetilde{U}}(M(\psi),M)\cong \Hom_{\widetilde{U}^{\geq 
0}}(V^{\psi},M)$, so there is a unique $\widetilde{U}$-module 
homomorphism $F:M(\psi)\rightarrow M$ with $F(v_{\psi})=m$, namely 
$F(u\otimes v)=u.m$ for all $u\in \widetilde{U}$.  \m For the final 
assertion, note that $\widetilde{U}^0$ acts semisimply on any 
finite-dimensional simple module $M$, and by (1.6) and Proposition 
1.7, there is some nonzero vector $m\in M_{\psi}$ such that $e_i.m=0 \ 
(1\leq i <n)$.  By the first part, $M$ is a quotient of $M(\psi)$, and 
so $M\cong L(\psi)$, as $L(\psi)$ is the unique simple quotient of 
$M(\psi)$.  \qed\enddemo \b As a special case, we will consider the 
modules $L(\lambda)=L(\hat{\lambda})$ where $\lambda \in \Lambda$.  
Let $\Lambda^+\subset \Lambda$ be the subset of {\it dominant} 
weights, that is $$\Lambda^+=\{\lambda\in\Lambda \mid 
\langle\alpha_i,\lambda\rangle \geq 0
\text{ for } 1\leq i <n\}.$$
We will show that if $L(\lambda)$ is finite-dimensional, then $\lambda\in
\Lambda^+$.  
This requires an identity for commuting $e_i$ past powers of $f_i$.  
For $k\geq 1$, let
$$[k]=\frac{r^k-s^k}{r-s}.  \tag 2.2$$
\b
\proclaim{Lemma 2.3} If $k\geq 1$, then $$\aligned e_if_i^k & 
=f_i^ke_i + [k]f_i^{k-1}\frac{r^{1-k}\omega_i - s^{1-k}\omega_i'} 
{r-s} \\
e_i^kf_i& =f_ie_i^k + [k]e_i^{k-1}\frac{s^{1-k}\omega_i - r^{1-k}\omega_i'}
{r-s}. \endaligned
$$ \endproclaim \m \demo{Proof}  For $k=1$, the above equations are just 
one of the defining relations of $U$.  Assume that
$k>1$ and $$
e_if_i^{k-1} = f_i^{k-1}e_i + [k-1]f_i^{k-2}\frac{r^{2-k}\omega_i-
s^{2-k}\omega_i'}{r-s}.$$
Then
$$
\aligned e_if_i^k & = \left(f_i^{k-1}e_i + [k-1]f_i^{k-2}\frac{r^{2-k}\omega_i
-s^{2-k}\omega_i'}{r-s}\right)f_i\\
&= f^{k-1}_i \left(f_ie_i + \frac{\omega_i-\omega_i'}{r-s}\right) + [k-1]
f_i^{k-1}\left(\frac{r^{1-k}s\omega_i - rs^{1-k}\omega_i'}{r-s}\right)\\
&= f_i^ke_i + \frac{f_i^{k-1}}{r-s}((1+[k-1]r^{1-k}s)\omega_i-(1+[k-1]rs^{1-k})
\omega_i')\\
&= f_i^ke_i + \frac{f_i^{k-1}}{r-s}([k]r^{1-k}\omega_i - [k]s^{1-k}\omega_i'). 
\endaligned $$

\n The argument for the second equation can be done similarly.  \qed 
\enddemo \b \proclaim{Lemma 2.4} Assume $rs^{-1}$ is not a root of unity.
Let $M$ be a nonzero 
finite-dimensional $\widetilde{U}$-module on which $\widetilde{U}^0$ 
acts semisimply, and $\l\in\Lambda$.  Suppose there is some nonzero 
vector $v\in M_{\l}$ with $e_i .  v=0$ for all $i \ (1\leq i <n)$.  
Then $\l\in\Lambda^+$.  \endproclaim \m \demo{Proof} Proposition 1.7 
implies that for any given value of $i$ there is some $k\geq 0$ such 
that $f_i^{k+1} .  v=0$ and $f_i^{k} .  v\neq 0$.  Applying $e_i$ to 
$f_i^{k+1} .  v=0$ and using Lemma 2.3 and the fact that $e_i .  v=0$, 
we have $$ 0 = [k+1] f_i^{k}\frac{r^{-k}\omega_i-s^{-k} 
\omega_i'}{r-s} .  v = 
\frac{[k+1]}{r-s}(r^{-k}\hat\lambda(\omega_i)-s^{-k}\hat\lambda(\omega_i'))
f_i^{k}.v.$$ 
Now $[k+1]/(r-s)\neq 0$ as $rs^{-1}$ is not a root of unity.  Therefore, since
$f^{k}.v\neq 0$,  
$$r^{-k}\hat\lambda(\omega_i)=s^{-k}\hat{\l}(\omega_i').$$
Equivalently, $$\aligned 
r^{-k}r^{\langle\epsilon_i,\l\rangle}s^{\langle\epsilon_{i+1},\l\rangle} 
&= 
s^{-k}r^{\langle\epsilon_{i+1},\l\rangle}s^{\langle\epsilon_{i},\l\rangle}, 
\\
\text{or} \ \ \ r^{-k + \langle\alpha_i,\l\rangle} &=
s^{-k + \langle\alpha_i,\l\rangle}. \endaligned $$
Again, because $rs^{-1}$ is not a root of unity, this forces $\langle\alpha_i,
\l\rangle = k \geq 0$.  Therefore $\l\in\Lambda^+$.  \qed \enddemo \b 

\proclaim {Corollary 2.5} When $rs^{-1}$ is not a root of unity, any 
finite-dimensional simple $\widetilde{U}$-module with weights in 
$\Lambda$ is isomorphic to $L(\lambda)$ for some $\lambda\in 
\Lambda^+$.  \endproclaim \b

 We will show next that 
all modules $L(\lambda)$ with $\l \in \Lambda^+$  are indeed 
finite-dimensional, and that all other finite-dimensional simple 
$\widetilde{U}$-modules are shifts of these by one-dimensional modules.  
In doing this, it helps to consider first the special case of simple 
$U_{r,s}({\fsl})$-modules.

\b \subhead Highest weight modules 
for $U = U_{r,s}(\fsl)$ \endsubhead 
\m 

For simplicity we drop the subscripts and just write $e,f,\w,\w'$ for the 
generators of $U= U_{r,s}(\fsl)$.  Any homomorphism $\phi: U^0 
\rightarrow \K$ is determined by its values on $\w$ and $\w'$.  By 
abuse of notation, we adopt the shorthand $\phi = \phi(\w)$ and $\phi' = 
\phi(\w')$.  \m Corresponding to each such $\phi$, there is a Verma
module $M(\phi) = U \ot_{U^{\geq 0}}\K v$ with basis $v_j = f^j \ot 
v \ \ (0 \leq j < \infty)$ such that the $U$-action is given by:
 
$$\aligned  f.v_j & = v_{j+1} \\
e.v_j & = [j]\frac{\phi r^{-j+1} -\phi's^{-j+1}}{r-s} v_{j-1} \qquad 
\qquad (v_{-1} := 0) \\
\w.v_j & = \phi r^{-j\la \e_1,\a_1\ra}s^{-j\la\e_2,\a_1\ra} v_j = \phi 
r^{-j}s^j v_j \\
\w'.v_j & = \phi' r^{-j\la\e_2,\a_1\ra}s^{-j\la\e_1,\a_1\ra} v_j = 
\phi' r^{j}s^{-j} v_j. \endaligned \tag 2.6$$

\n Note that $M(\phi)$ is a simple $U$-module if and only if $\displaystyle 
{[j]\frac{\phi r^{-j+1} -\phi's^{-j+1}}{r-s}} \neq 0$ for any $j \geq 
1$.  \m Suppose $\displaystyle {[\ell+1]\frac{\phi r^{-\ell} 
-\phi's^{-\ell}}{r-s} = 0}$ for some $\ell \geq 0$.  Then either 
$r^{\ell+1} = s^{\ell+1}$, which implies $rs^{-1}$ is a root of unity, 
or $\phi' = \phi r^{-\ell} s^\ell$.  Assuming that $rs^{-1}$ is not a 
root of unity and $\phi' = \phi r^{-\ell}s^\ell $, we see that the 
elements $v_i, \ i \geq \ell+1$, span a maximal submodule.  The 
quotient is the $(\ell+1)$-dimensional simple module $L(\phi)$, which 
we can suppose is spanned by $v_0,v_1, \dots, v_\ell$ and has 
$U$-action given by

$$\aligned f.v_j & = v_{j+1}, \qquad \quad (v_{\ell+1} = 0) \\
e.v_j & =  \phi r^{-\ell}[j][\ell+1-j] v_{j-1} \qquad \qquad
(v_{-1} = 0) \\
\w.v_j & = 
\phi r^{-j}s^j v_j \\
\w'.v_j & = \phi r^{-\ell +j}s^{\ell-j} v_j.  \endaligned \tag 2.7$$

\m When $M(\phi)$ is not simple and $rs^{-1}$ is
not a root of unity, $j = \ell+1$ is the unique value such that 
$\displaystyle {[j]\frac{\phi r^{-j+1} -\phi's^{-j+1}}{r-s} = 0}$. 
In this case, $M(\phi)$ has a unique proper submodule, namely the maximal 
submodule generated by $v_{\ell +1}$ as above.
 
\m
 We now have the following classification of simple modules for 
$U_{r,s}({\frak{sl}}_2)$.  

\b \proclaim {Proposition 2.8} 
\roster
\item"{(i)}"  Assume $U = U_{r,s}(\fsl)$, where $rs^{-1}$ is not a root of unity.  
Let $\phi: U^0 \rightarrow \K$ be an algebra homomorphism such that 
$\phi(\w') = \phi(\w)r^{-\ell}s^\ell$ for some $\ell \geq 0$.  Then 
there is an $(\ell+1)$-dimensional simple $U$-module $L(\phi)$ spanned 
by vectors $v_0,v_1,\dots, v_\ell$ and having $U$-action given by 
(2.7).  Any $(\ell+1)$-dimensional simple $U$-module is isomorphic to 
some such $L(\phi)$.
\m \item"{(ii)}" If $\nu = \nu_1\e_1 + \nu_2 \e_2 \in \Lambda^+$, then 
$\nu_1-\nu_2 = \ell$ for some $\ell \in \Z_{\geq 0}$, and $\nu(\w') = 
r^{\nu_2}s^{\nu_1} = r^{\nu_1-\ell}s^{\nu_2+\ell} = 
\nu(\w)r^{-\ell}s^\ell$ in this case.  Thus, the module $L(\nu)$ is 
$(\ell+1)$-dimensional and has $U$-action given by (2.7) with $\phi = 
r^{\nu_1}s^{\nu_2} = r^{\nu_1}s^{\nu_1-\ell}$.  \endroster \endproclaim

\b \subhead Finite-dimensionality of $L(\l)$ 
for $\l \in \Lambda^+$ \endsubhead

\m
We show below that the simple $\widetilde U$-modules $L(\lambda)$ with 
$\lambda\in \Lambda^+$ are finite-dimensional.  For this it suffices 
to prove that $M(\l)$ has a $\widetilde{U}$-submodule of finite 
codimension, as $L(\l)$ is the quotient of $M(\l)$ by its unique 
maximal submodule.  \m As $\l$ is dominant, $k_i=\langle 
\alpha_i,\l\rangle$ for $i=1,\ldots,n-1$, are nonnegative integers.  
Define a $\widetilde{U}$-submodule $M'(\l)$ of $M(\l)$ by
$$M'(\l)=\sum_{i=1}^{n-1} \widetilde{U}f_i^{k_i+1} .  v_{\l}.  \tag 2.9$$
Our goal is to prove that the module $L'(\lambda)=M(\lambda)/M'(\lambda)$ 
is nonzero and finite-dimensional.
\m
 By Lemma 2.3 we have $$\aligned e_if_i^{k_i+1}.v_{\l} &= 
[k_i+1]f_i^{k_i}\frac{r^{-k_i}\omega_i - s^{-k_i}\omega_i'}{r-s}. 
v_{\l}\\
&= [k_i+1] f_i^{k_i} \frac{r^{-\langle\alpha_i,\l\rangle}
 r^{\langle\epsilon_i,\l\rangle} s^{\langle\epsilon_{i+1},\l\rangle} - 
 s^{-\langle\alpha_i,\l\rangle}r^{\langle\epsilon_{i+1},\l
  \rangle}s^{\langle\epsilon_i,\l\rangle}}{r-s}.v_{\l}\\
&= [k_i+1]f_i^{k_i}\frac{r^{\langle\epsilon_{i+1},\l\rangle}s^{\langle
 \epsilon_{i+1},\l\rangle} - s^{\langle\epsilon_{i+1},\l\rangle}
  r^{\langle\epsilon_{i+1},\l\rangle}}{r-s}.v_{\l} = 0. \endaligned $$
If $j\neq i$, $e_jf_i^{k_i+1}.v_{\l}=f_i^{k_i+1}e_j.v_{\l}=0$ by the 
defining relations.  Consequently, by Theorem 2.1, 
$\widetilde{U}f_i^{k_i+1}.v_{\l}$ is a homomorphic image of $M(\l 
-(k_i+1)\alpha_i)$, and so all its weights are less than or equal to 
$\lambda - (k_i+1)\alpha_i$.  This implies that $v_{\l}\not\in 
M'(\l)$, hence $L'(\l) \neq 0$.
 
\b

\proclaim{Lemma 2.10} The elements $e_j,f_j \ (1\leq j<n)$ act locally nilpotently
on $L'(\l)$. \endproclaim

\m

\demo{Proof} As the Verma module $M(\l)$ is spanned over $\K$ by all elements 
$x_1\cdots x_t.v_{\l}$ where $x_1,\ldots,x_t\in \{f_1,\ldots,f_{n-1}\}, 
\ t\in \Z_{\geq 0}$, it is enough to argue by induction on $t$ that a 
sufficiently high power of $e_j$ (resp., $f_j$) takes such an element 
to $M'(\l)$.  If $t=0$, then $e_j.v_{\l}=0\in M'(\l)$, and 
$f_j^{k_j+1}.v_{\l}\in M'(\l)$ by construction.  Now assume that there 
are positive integers $N_j$ such that $$e_j^{N_j}x_2\cdots x_t 
.v_{\l}\in M'(\l) \ \text{ and } \
f_j^{N_j}x_2\cdots x_t.v_{\l}\in M'(\l).$$
Suppose that $x_1=f_i$.  If $j\neq i$, then  
$$e_j^{N_j}x_1\cdots x_t.v_{\l} =f_ie_j^{N_j}x_2\cdots x_t.v_{\l}\in M'(\l).$$
Otherwise by Lemma 2.3,
$$e_i^{N_i+1}x_1\cdots x_t.v_{\l}=f_ie_i^{N_i+1}x_2\cdots x_t.v_{\l}
+[N_i+1] e_i^{N_i}\frac{s^{-N_i}\omega_i - r^{-N_i}\omega_i'}{r-s}
x_2\cdots x_t.v_{\l}.$$
Applying relation (R2') and the induction hypothesis, we see that these
terms are both in $M'(\l)$.
\m
Now $f_i^{N_i -1} x_1\cdots x_t.v_{\l}=f_i^{N_i}x_2\cdots x_t.v_{\l} 
\in M'(\l)$, and if $|i-j|>1$, we also have $f_j^{N_j}x_1\cdots 
x_t.v_{\l} =f_if_j^{N_j}x_2\cdots x_t.v_{\l} \in M'(\l)$.  Finally, we 
need to show that if $|i-j|=1$, then $f_j^{N_j+1}x_1\cdots x_t.v_{\l}\in M'(\l)$.  This will follow from the induction hypothesis 
once we know that $f_j^{N_j+1}f_i\in \K f_jf_if_j^{N_j}+\K 
f_if_j^{N_j+1}$.  \m We argue by induction on $m\geq 1$ that
$$f_j^{m+1}f_i\in \K f_jf_if_j^m + \K f_if_j^{m+1}.$$
Indeed if $m=1$, this follows from relation (R7), but if $m>1$, then by 
induction and (R7), $$f_j^{m+1}f_i\in f_j(\K f_jf_if_j^{m-1} + \K 
f_if_j^m)\subseteq \K f_jf_if_j^m + \K f_if_j^{m+1}.  \qed $$ \enddemo

\m

\proclaim{Lemma 2.11} Assume $rs^{-1}$ is not a root of unity, and
let $V$ be a 
module for $U = U_{r,s}({\frak{sl}}_2)$ on which $U^0$ acts 
semisimply.  Assume $V=\oplus_{j\in \Z_{\geq 0}} V_{\l - j\alpha}$ for 
some weight $\l \in \Lambda$, each weight space of $V$ is 
finite-dimensional, and $e$ and $f$ act locally nilpotently on $V$.  
Then $V$ is finite-dimensional, and the weights of $V$ are preserved 
under the simple reflection taking $\alpha$ to $-\alpha$.  
\endproclaim

\m

\demo{Proof}  Let $\mu=\mu_1\epsilon_1 + \mu_2\epsilon_2$ be a weight of $V$,
and $v\in V_{\mu}\setminus \{0\}$.  As $e$ acts locally nilpotently on $V$,
there is a nonnegative integer $k$ such that $e^{k+1}.v=0$ and $e^k.v\neq 0$.
By Theorem 2.1, $Ue^k.v$ is a homomorphic image of $M(\mu+k\alpha)$.  
But since $f$ acts locally nilpotently on $Ue^k.v$, this image cannot be 
isomorphic to $M(\mu+k\alpha)$.  Thus  because $M(\mu+k\alpha)$ has
a unique proper submodule, 
$U e^k.v \cong L(\mu +k\alpha)$, and so it is finite-dimensional.  
Corollary 2.5 implies that $\mu+k\alpha$ is dominant.  As there are only 
finitely many dominant weights less than or equal to the given weight 
$\l$, and each weight space is finite-dimensional, it must be that $V$ 
itself is finite-dimensional.

\m In particular, $V$ has a composition series with factors isomorphic to 
$L(\nu)$ for some $\nu \in \Lambda^+$.  Any weight $\mu$ of $V$ is a 
weight of some such $L(\nu)$ with $\nu = \nu_1\e_1 + \nu_2 \e_2 \in 
\Lambda^+$.  By (ii) of Proposition 2.8, $L(\nu)$ has weights $\nu, 
\nu-\a, \dots, \nu-\ell \a$ where $\ell = \nu_1-\nu_2$.  Thus, $\mu = 
\nu-j\a$ for some $j \in \{0,1,\dots,\ell\}$.  But then $\mu-\la 
\mu,\a\ra \a = \nu-(\ell-j)\a$ is a weight of $L(\nu)$ since $\ell-j 
\in \{0,1,\dots,\ell\}$, hence it is a weight of $V$.  Thus, the 
weights of $V$ are preserved under the simple reflection taking 
$\alpha$ to $-\alpha$.  \qed\enddemo

\b

\proclaim{Lemma 2.12} Assume that $rs^{-1}$ is not a root of unity, and
let $\l\in \Lambda^+$.  Then $L(\l)$ is finite-dimensional.
\endproclaim

\m

\demo{Proof} This follows once we show that $L'(\l) = M(\l)/M'(\l)$,
where $M'(\l)$ is as in (2.9), is finite-dimensional.  We will prove that 
the set of weights of $L'(\l)$ is preserved under the action of the 
symmetric group $S_n$ (the Weyl group of $\fgl$) on $\Lambda$ which
is generated by the simple reflections $s_i: \mu \rightarrow \mu - \la 
\mu,\a_i \ra \a_i \ (1\leq i <n)$.  Each $S_n$-orbit contains a 
dominant weight, and there are only finitely many dominant weights 
less than or equal to $\l$.  As the weights in $M(\l)$ are all less 
than or equal to $\l$, and the weight spaces are finite-dimensional, 
the same is true of $L'(\l)$.  Therefore $L'(\l)$ is 
finite-dimensional.  \m To see that $s_i$ preserves the set of weights 
of $L'(\l)$, let $\mu=\mu_1\epsilon_1+\cdots +\mu_n\epsilon_n$ be a 
weight of $L'(\l)$.  Consider $L'(\l)$ as a module for the copy $U_i$ 
of $U_{r,s}({\frak {sl}}_2)$ generated by $e_i,f_i, \omega_i, 
\omega_i'$, and let $L'_i(\mu)$ be the $U_i$-submodule of $L'(\l)$ 
generated by $L'(\l)_{\mu}$.  As all weights of $L'(\l)$ are less than 
or equal to $\l$, we have
$$L_i'(\mu) = \bigoplus_{j\in\Z_{\geq 0}} L_i'(\mu)_{\lambda'-j\alpha_i}$$
for some weight $\lambda' \leq \lambda$.  By 
Lemmas 2.10 and 2.11, the simple reflection $s_i$ preserves the 
weights of $L'_i(\mu)$, so in particular, $s_i(\mu)$ is also a weight 
of $L'(\l)$.  \qed\enddemo

\b
\n {\bf Remark 2.13.} It will follow from Lemma 3.7 in the next 
section that $L(\l)\cong L'(\l)$, since $L(\l)$ is the unique simple 
quotient of $M(\l)$, $L'(\l)$ is a finite-dimensional quotient of 
$M(\l)$, and by that lemma, every finite-dimensional quotient is 
simple.

\b

\proclaim{Corollary 2.14}  Assume that $rs^{-1}$ is not a root of unity.
The finite-dimensional simple $\widetilde{U}$-modules having weights in
$\Lambda$ are precisely the modules $L(\l)$ where $\l\in \Lambda^+$.
Moreover, $L(\l)\cong L(\mu)$ if and only if $\l=\mu$.
\endproclaim

\m

\demo{Proof}
The first statement is a consequence of Corollary 2.5 and 
Lemma 
2.12.  Assume there is an isomorphism of $\widetilde{U}$-modules from 
$L(\l)$ to $L(\mu)$.  The highest weight vector of $L(\l)$ must be 
sent to a weight vector of $L(\mu)$, so $\l\leq \mu$.  As a similar 
argument shows that $\mu\leq \l$, we have $\l=\mu$.  \qed\enddemo

\b \subhead Shifts by one-dimensional modules 
\endsubhead \m Suppose now that we have a one-dimensional module $L$ 
for $\widetilde{U} = U_{r,s}(\fgl)$.  Then by Theorem 2.1, $L=L(\chi)$ 
for some algebra homomorphism $\chi:\widetilde{U}^0 \rightarrow \K$, 
with the elements $e_i, f_i \ (1\leq i <n)$ acting as multiplication 
by 0.  Relation (R4) yields
$$
\chi(\w_i) =  \chi(a_ib_{i+1})=\chi(a_{i+1}b_i) = \chi(\w_i') \qquad  
(1\leq i<n). \tag 2.15$$
Conversely, if an algebra homomorphism $\chi$ satisfies this equation, then 
$L(\chi)$ is one-dimensional by relation (R4).  We will write 
$L_{\chi}=L(\chi)$ to emphasize that the module is one-dimensional. 
\b \proclaim 
{Proposition 2.16} Assume $rs^{-1}$ is not a root of unity and 
$L(\psi)$ is the finite-dimensional simple module for $\widetilde U = 
U_{r,s}(\fgl)$ with highest weight $\psi$.  Then there exists a 
homomorphism $\chi: \widetilde U^0 \rightarrow \K$ such that (2.15) 
holds and an element $\l \in \Lambda^{+}$ so that $\psi = \chi \cdot \hat 
\l$.  Thus, the weights of $L(\psi)$ belong to $\chi \cdot \hat 
\Lambda$.  \endproclaim \m \demo {Proof} When $L(\psi)$ is viewed as a 
module for the copy $U_i$ of $U_{r,s}(\fsl)$ generated by $e_i, f_i, 
\w_i, \w_i'$, it has a composition series whose factors are simple 
$U_i$-modules as described by Proposition 2.8.  As the highest weight 
vector of $L(\psi)$ gives a highest weight vector of some composition 
factor, there is a weight $\phi_i$ of $U_i$ and a nonnegative integer 
$\ell_i$ so that $\psi(\w_i) = \phi_i(\w_i)$ and $\psi(\w_i') = 
\phi_i(\w_i') = \phi_i(\w_i)r^{-\ell_i}s^{\ell_i} = 
\psi(\w_i)r^{-\ell_i}s^{\ell_i}$.
\m Set $\ell_n = 0$ and define $\l_i = \ell_i + \dots + \ell_n$ for $i 
=1,\dots, n$.  Let $\l = \sum_{i=1}^n \l_i \e_i$, which belongs to 
$\Lambda^+$.  Now we define $\chi: \widetilde U^0 \rightarrow \K$ by the 
formulas

$$\aligned 
\chi(a_i) & = \psi(a_i)r^{-\la \e_i,\l\ra} = \psi(a_i)r^{- ( \ell_i+\cdots
+\ell_n )} \\
\chi(b_i) & = \psi(b_i)s^{-\la \e_i,\l\ra} =
\psi(b_i)s^{-( \ell_i+\cdots+\ell_n)}.\endaligned$$

\n Then it follows that 

$$\aligned \chi(\w_i') & = \chi(a_{i+1}b_{i}) = \psi(\w_i')r^{-(\ell_{i+1}+\cdots+\ell_n)}
s^{-(\ell_i+\cdots+\ell_n)} \\
& = \psi(\w_i)r^{-\ell_i}s^{\ell_i}r^{-(\ell_{i+1}+\cdots+\ell_n)}
s^{-(\ell_i+\cdots+\ell_n)} \\
& = \chi(a_ib_{i+1}) = \chi(\w_i) \endaligned$$

\n for $i =1,\dots, n-1$, and $\psi = \chi \cdot \hat \l$ as desired.  \qed \enddemo
\b
\n {\bf Remark 2.17.} If $M$ is any finite-dimensional module,
then $M = \bigoplus_{i=1}^m \bigoplus_{\l \in \Lambda} M_{\psi_i \cdot \hat \l}$
for some weights $\psi_i$ such that $\psi_i \cdot \hat \Lambda$\ $(1 \leq i \leq m)$
are distinct cosets in $\Hom(\widetilde U^0, \K)/\hat \Lambda$ (viewed as a $\Z$-module
under the action $k \cdot \psi = \psi^k$).   Then $M_i :=  \bigoplus_{\l \in \Lambda}
M_{\psi_i \cdot \hat \l}$ is a submodule, and $M = \bigoplus_{i=1}^m M_i$. Therefore,
if $M$ is an indecomposable $\widetilde U$-module, $M = \bigoplus_{\l \in \Lambda}
M_{\psi \cdot \hat \l}$  for some $\psi \in \Hom(\widetilde U^0, \K)$. 
A simple submodule $S$ of $M$ has weights in $\psi \cdot \hat \Lambda$. By
replacing $\psi$ with the homomorphism $\chi$ for $S$ given by Proposition
2.16, we may assume that for any indecomposable module $M$, there
is a $\chi$ satisfying (2.15) so that $M = \bigoplus_{\l \in \Lambda} M_{\chi \cdot \hat \lambda}$.

\b
\proclaim{Lemma 2.18} Let $\chi:\widetilde{U}^0\rightarrow \K$ be an algebra
homomorphism with $\chi(\w_i) =\chi(\w_i') \ \ (1\leq
i<n)$. Let $M$ be a finite-dimensional $\widetilde{U}$-module whose weights are all
in $\chi\cdot\hat{\Lambda}$.  If $\widetilde{U}^0$ acts semisimply
on $M$, then
$$
  M\cong L_{\chi}\otimes N
$$
for some $\widetilde{U}$-module $N$ whose weights are all in 
$\Lambda$.
\endproclaim
\m \demo{Proof} Let $\chi^{-1}:\widetilde{U}^0\rightarrow \K$ be the 
algebra homomorphism defined by 
$\chi^{-1}(a_i)=\chi(a_i^{-1})=(\chi(a_i))^{-1}$ and 
$\chi^{-1}(b_i)=\chi(b_i^{-1})=(\chi(b_i))^{-1}$ for $1\leq i \leq n$.  
Note that $L_{\chi}\otimes L_{\chi^{-1}}$ is isomorphic to the trivial 
module $L_{\varepsilon}$ corresponding to the counit.  Let
$$N=L_{\chi^{-1}}\otimes M.$$
Then $M\cong L_{\chi}\otimes N$ as $L_{\varepsilon}$ is a multiplicative 
identity (up to isomorphism) for $\widetilde{U}$-modules.  
The weights of $N$ are all in $\chi^{-1}\cdot \chi\cdot\hat{\Lambda}=
\hat{\Lambda}$. \qed \enddemo
\b

We now have a classification of finite-dimensional simple 
$\widetilde{U}$-modules.

\b \proclaim{Theorem 2.19} Assume $rs^{-1}$ is not a root of unity.
The finite-dimensional simple $\widetilde{U}$-modules are precisely the
modules
$$ L_{\chi}\otimes L(\l),$$
where $\chi: \widetilde{U}^0\rightarrow \K$ is an algebra homomorphism with 
$\chi(\omega_i) =\chi(\omega_i') \ (1\leq i <n)$, and 
$\l\in\Lambda^+$.  \endproclaim \b \demo{Proof} Let $M$ be a 
finite-dimensional simple $\widetilde{U}$-module.  By Theorem 2.1, 
Proposition 2.16, and Lemma 2.18, $M\cong L_{\chi}\otimes N$ for some 
$\chi$ satisfying (2.15) and some simple module $N$ with weights in 
$\Lambda$.  By Corollary 2.5, $N\cong L(\lambda)$ for some 
$\lambda\in\Lambda^+$.  Conversely, any $\widetilde{U}$-module of this 
form is finite-dimensional by Lemma 2.12 and simple
by its construction.  \qed \enddemo \b \n {\bf Remark 2.20.} If $r=q$ and 
$s=q^{-1}$ for some $q\in\K$, the classification of finite-dimensional 
simple $U_q({\frak {sl}}_n)$-modules is a consequence of Theorem 2.19 
applied to $U_{q,q^{-1}}({\frak {sl}}_n)$: \ The simple $U_q({\frak 
{sl}}_n)$-modules are precisely those simple $U_{q,q^{-1}}({\frak 
{sl}}_n)$-modules on which $\omega_i'$ acts as $\omega_i^{-1}$, so 
that
$$\chi(\omega_i)=\chi(\omega_i')=\chi(\omega_i^{-1}).$$
This implies $\chi(\omega_i)=\pm 1 \ (1\leq i <n)$.  Each choice of algebra
homomorphism $\chi: U^0\rightarrow \K$ with $\chi(\omega_i)=\chi(\omega'_i)
=\pm 1$
yields a one-dimensional $U_{q,q^{-1}}({\frak {sl}}_n)$-module $L_{\chi}$,
and so the simple $U_q({\frak{sl}}_n)$-modules are the $L_{\chi}\otimes
L(\lambda)$ with $\lambda\in\Lambda^+$ and $\chi$ as above. 
(Compare with [Ja, \S 5.2, Convention 5.4, and Thm. 5.10].) \b \n{\bf 
Remark 2.21.} We can interpret Proposition 2.8 in light of Theorem 
2.19: \ Let $L(\phi)$ be the simple $U_{r,s}({\frak{sl}}_2)$-module 
described in the proposition.  Let $\lambda=\ell\epsilon_1\in 
\Lambda^+$ and define $\chi: U^0\rightarrow\K$ by 
$\chi(\omega)=\phi(\omega) r^{-\ell}$, 
$\chi(\omega')=\phi(\omega')s^{-\ell}=\phi(\omega)r^{-\ell}s^{\ell}s^{-\ell} 
=\chi(\omega)$.  Then $\phi=\chi\cdot \hat{\lambda}$ and $L(\phi)\cong 
L_{\chi}\otimes L(\lambda)$.

\b \m

\head \S 3. Complete reducibility \endhead \m
In this section
we will establish complete reducibility of all 
finite-dimensional $\widetilde U$-modules on which $\widetilde U^{0}$ 
acts semisimply.  However, it is helpful to
 work in a more general 
context.  \m Let $\Cal O$ denote the category of modules $M$ for 
$\widetilde U = U_{r,s}(\fgl)$ which satisfy the conditions: \m \roster 
\item"{($\Cal O$1)}" $\widetilde U^0$ acts semisimply on $M$, and the 
set $\text{wt}(M)$ of weights of $M$ belongs to $\Lambda$: \quad $M = 
\bigoplus_{\l \in \text{wt}(M)} M_{\l}$, where $M_{\l}= \{m \in M \mid 
a_i.m = r^{\la \e_i,\l\ra}, \ \ b_i.m = s^{\la \e_i,\l \ra}$ for all 
$i\}$; \item"{($\Cal O$2)}" $\dim _{\K} M_{\l} < \infty$ for all $\l \in 
\text{wt}(M)$; \item"{($\Cal O$3)}" $\text{wt}(M) \subseteq 
\bigcup_{\mu \in F} (\mu - Q^+)$ for some finite set $F \subset 
\Lambda$.  \endroster \m \n The morphisms in $\Cal O$ are $\widetilde 
U$-module homomorphisms.  \m  
All finite-dimensional $\widetilde 
U$-modules which satisfy (1) belong to category $\Cal O$, as do all 
highest weight modules with weights in $\Lambda$ such as the Verma 
modules $M(\l)$.

\m
 We recall the definition of the quantum Casimir operator [BW, 
Sec.  4].  It is a consequence of (R2) and (R3) that the subalgebra 
$U^+$ of $\widetilde{U}$ (or of $U=U_{r,s}({\frak{sl}}_n)$) generated 
by 1 and $e_i \ (1\leq i <n)$ has the decomposition 
${U}^+=\oplus_{\zeta\in Q^+} {U}^+_{\zeta}$ where 
$${U}^+_{\zeta}=\{z\in U^+\mid a_iz=r^{\langle\epsilon_i , \zeta 
\rangle }za_i, \ b_iz=s^{\langle\epsilon_i , \zeta\rangle }zb_i \
(1\leq i < n)\}.$$
The weight space ${U}^+_{\zeta}$ is spanned by all the monomials
$e_{i_1}\cdots e_{i_{\ell}}$ such that $\alpha_{i_1} +\cdots +
\alpha_{i_{\ell}} =\zeta$.  Similarly, the subalgebra ${U}^-$
generated by 1 and the $f_i$ has the decomposition ${U}^-=
\oplus_{\zeta\in Q^+}{U}^-_{-\zeta}$. The spaces ${U}^+_{\zeta}$ 
and ${U}^-_{-\zeta}$ are nondegenerately paired by
the Hopf pairing defined by 
$$\aligned
(f_i,e_j) &= \frac{\delta_{i,j}}{s-r}\\
(\omega_i',\omega_j) &= r^{\langle \epsilon_j,\alpha_i\rangle}
s^{\langle\epsilon_{j+1},\alpha_i\rangle}\\
(b_n,a_n) & = 1, \ \ \ (b_n,\omega_j)=s^{-\langle\epsilon_n,\alpha_j\rangle},
\ \ \ (\omega_i',a_n)=r^{\langle \epsilon_n,\alpha_i\rangle }.
\endaligned \tag 3.1 $$
(See [BW, Sec. 2].)  The Hopf algebras $\widetilde{U}$ and $U$ are 
Drinfel'd doubles of certain Hopf subalgebras with respect to this 
pairing [BW, Thm.  2.7].  Let $d_{\zeta}=\dim_{\K}U^+_{\zeta}$.  Assume 
$\{u^{\zeta}_k\}^{d_{\zeta}}_{k=1}$ is a basis for $U^+_{\zeta}$, and 
$\{v^{\zeta}_k\}^{d_{\zeta}}_{k=1}$ is the dual basis for $U^-_{-\zeta}$
with respect to the pairing.
\m
 Now let $$\Omega=\sum_{\zeta\in Q^+}
\sum_{k=1}^{d_{\zeta}}S(v^{\zeta}_k)u^{\zeta}_k, \tag 3.2$$
where $S$ denotes the antipode. 
All but finitely many terms in this sum will
act as multiplication by 0 on any  
weight space $M_\l$ of $M \in \Cal O$.  Therefore $\Omega$ is a well-defined 
operator on such $M$.

\m
 The second part of the Casimir operator involves a function 
$g:\Lambda \rightarrow \K^\#$ defined as follows. 
If $\rho$ denotes the half sum of the positive roots,
then $2\rho  
 = \sum_{j=1}^n (n+1-2j)\e_j \in \Lambda$.  For 
$\lambda \in \Lambda$, set

$$g(\lambda) = (rs^{-1})^{\half \la \lambda + 2 \rho, \lambda 
\ra}. \tag 3.3$$

\n When $M$ is a $\widetilde{U}$-module in $\Cal O$, we define the linear 
operator $\Xi :M\rightarrow M$ by
$$ \Xi (m) = g(\l)m$$
for all $m\in M_{\l}, \ \l\in\Lambda$.  Then we have the following 
result from [BW].  \b \proclaim{Proposition 3.4} [BW, Thm. 4.20] \ The 
operator $\Omega\Xi:M\rightarrow M$ commutes with the action of 
$\widetilde{U}$ on any $\widetilde U$-module $M \in \Cal O$.  
\endproclaim \b We require the next lemma in order to prove complete 
reducibility.

  \b \proclaim {Lemma 3.5} Assume 
$rs^{-1}$ is not a root of unity, and let $\lambda,\mu \in 
\Lambda^{+}$.  If $\lambda \geq \mu$ and $g(\lambda) = g(\mu)$, then 
$\lambda = \mu$.  \endproclaim

 \m \demo {Proof} \ Because $\lambda \geq \mu$, we may suppose $\lambda = \mu + 
\beta$ where $\beta = \sum_{i=1}^{n-1}k_i \a_i$ and $k_{i} \in 
\Z_{\geq 0}$.  By assumption we have

$$(rs^{-1})^{{\half} \la \lambda + 2 \rho, \lambda \ra} = g(\lambda) = 
g(\mu) = (rs^{-1})^{{\half} \la \mu + 2 \rho, \mu
\ra},$$  

\n and as $rs^{-1}$ is not a root of unity, it must be that
$\la \lambda+2\rho,\lambda\ra = \la \mu+2\rho,\mu \ra$, or equivalently, 
$2 \la \mu+\rho, \beta \ra + \la \beta,\beta \ra = 0$.  Since $\mu \in 
\Lambda^{+}$, $\mu = \mu_1\e_1 + \mu_2 \e_2 + \cdots + \mu_n \e_n$ 
where $\mu_i \in \Z$ for all $i$ and $\mu_1 \geq \mu_2 \geq \cdots 
\geq \mu_n$.  Then

$$
\eqalign {0 & = \la 2\mu+2\rho, \beta \ra + \la \beta,\beta \ra \cr & = 
\sum_{i=1}^{n-1} k_i\Big(2\mu_i +(n+1-2i) 
-2\mu_{i+1}-(n+1-2(i+1))\Big) + \sum_{i=1}^{n} (k_i-k_{i-1})^2 \cr & 
\hskip 4.1 truein (k_0 = 0 = k_n) \cr & = \sum_{i=1}^{n-1} 
2k_i\big(\mu_i-\mu_{i+1}+1 \big) + \sum_{i=1}^{n}
(k_i-k_{i-1})^2.}$$

\n The only way this can happen is if $k_i = 0$ for all $i$ and
$\lambda = \mu$.   
 \qed \enddemo 
 \b \proclaim{Lemma 3.6} Assume that $rs^{-1}$ is not a root of unity.
\roster 
  \item"{(i)}" $\Omega \Xi$ acts as multiplication by the scalar $g(\l) 
 = (rs^{-1})^{\half \la \lambda + 2 \rho, \lambda \ra}$ on the Verma 
 module $M(\l)$ with $\l\in\Lambda$, hence on any submodule or quotient 
 of $M(\lambda)$.
 \item"{(ii)}" The eigenvalues of the operator $\Omega\Xi:M\rightarrow 
 M$ are integral powers of $(rs^{-1})^{\half}$ on any finite-dimensional 
 $M \in \Cal O$.  \endroster \endproclaim \m \demo{Proof} By its 
 construction, $\Omega \Xi$ acts by multiplication by $g(\l) = 
 (rs^{-1})^{\half \la \lambda + 2 \rho, \lambda \ra}$ on the maximal 
 vector $v_{\l}$ of $M(\l)$.  But since $M(\l) = \widetilde U . v_{\l}$ 
 and $\Omega \Xi$ commutes with $\widetilde U$ on modules in $\Cal O$, 
 $\Omega \Xi$ acts as multiplication by $(rs^{-1})^{\half \la \lambda + 
 2 \rho, \lambda \ra}$ on all of $M(\l)$.  \m If $M \in \Cal O$ is 
 finite-dimensional, it has a composition series.  Each factor is a 
 finite-dimensional simple $\widetilde{U}$-module with weights in 
 $\Lambda$, and in particular, is a quotient of $M(\l)$ for some $\l\in 
 \Lambda$.  On such a factor, $\Omega \Xi$ acts as multiplication by 
 $g(\l)$.  Therefore the action of $\Omega \Xi$ on $M$ may be expressed 
 by an upper triangular matrix with each diagonal entry equal to 
 $g(\l)$ for some $\l\in\Lambda$.  \qed\enddemo

 \b \proclaim{Lemma 3.7} Assume $rs^{-1}$ is not a root of unity.
Let $\l\in\Lambda$ and $M$ be a nonzero finite-dimensional 
quotient of the Verma module $M(\l)$.  Then $M$ is simple.  
\endproclaim \m \demo{Proof} First observe that by Lemma 2.4, 
$\l\in\Lambda^+$.  Assume $M'$ is a proper submodule of $M$.  As $M$ 
is generated by its one-dimensional subspace $M_{\l}$, we must have 
$M_{\l}'=0$.  Let $\mu\in\Lambda$ be maximal such that $M_{\mu}'\neq 
0$, and note that $\mu<\l$.  Let $m'$ be a nonzero vector of 
$M_{\mu}'$.  By maximality of $\mu$, we have $e_i.m'=0$ for all $i \ 
(1\leq i<n)$.  Letting $M''=U.m'$, a nonzero finite-dimensional 
quotient of $M(\mu)$, we see that $\mu\in\Lambda^+$ as well.  By Lemma 
3.6 (i), $\Omega\Xi$ acts as multiplication by $g(\l)$ on $M$, and by 
$g(\mu)$ on $M''$.  This forces $g(\l)=g(\mu)$, which contradicts 
Lemma 3.5 as $\mu<\l$.  \qed\enddemo \b \proclaim{Theorem 3.8} Assume 
$rs^{-1}$ is not a root of unity.  Let $M$ be a nonzero 
finite-dimensional $\widetilde{U}$-module on which $\widetilde{U}^0$ 
acts semisimply.  Then $M$ is completely reducible.  \endproclaim \m 
\demo{Proof} We will establish the result first in the case $M$ has 
weights in $\Lambda$.  Write $M$ as a direct sum of generalized 
eigenspaces for $\Omega\Xi$.  Note that by Proposition 3.4, this is a 
direct sum decomposition of $M$ as a $\widetilde{U}$-module.  
Therefore we may assume $M$ is itself a generalized eigenspace of 
$\Omega\Xi$, so that $(\Omega\Xi-(rs^{-1})^c)^d(M)=0$ for some $c\in 
\half\Z, \ d=\dim_{\K}M$, by Lemma 3.6 (ii).  \m Let $P=\{m\in M\mid 
e_i .  m=0 \ (1\leq i <n)\}$, and note that $P=\oplus 
_{\l\in\Lambda}P_{\l}$, $P_{\l}=P\cap M_{\l}$.  If $m\in 
P_{\l}-\{0\}$, the $\widetilde{U}$-submodule $\widetilde{U}.m$ of $M$ 
is a nonzero quotient of $M(\l)$ by Theorem 2.1.  By Lemma 3.7, each 
such $\widetilde{U}.m$ is a simple $\widetilde{U}$-module, and so the 
$\widetilde{U}$-submodule $M'$ of $M$ generated by $P$ is a sum of 
simple $\widetilde{U}$-modules.  That is, $M'$ is completely 
reducible.  Let $M''=M/M'$.  \m Assuming $M''\neq 0$, there is a 
weight $\mu\in\Lambda$ maximal such that $M_{\mu}''\neq 0$.  Let 
$m''\in M''_{\mu}-\{0\}$.  By maximality of $\mu$, we have $e_i.m''=0$ 
for all $i \ (1\leq i<n)$.  By Lemma 2.4, we have $\mu\in\Lambda^+$, 
and by Theorem 2.1 and Lemma 3.6, $\Omega\Xi$ acts as multiplication 
by $g(\mu)$ on the $U$-module $U.m''$ generated by $m''$.  This 
implies $g(\mu)= (rs^{-1})^c$.  \m Let $m\in M_{\mu}$ be a 
representative for $m''\in (M/M')_{\mu}$, and $M_1=\widetilde{U} .m$.  
Then the module $M_1$ is a direct sum of its intersections with the 
weight spaces of $M$, so there is a weight $\eta\in \Lambda$ maximal 
such that $M_1\cap M_{\eta}\neq 0$.  Let $m_1\in M_1\cap 
M_{\eta}-\{0\}$, so that $e_i.m_1=0$ for all $i \ (1\leq i<n)$.  Again 
applying Theorem 2.1 and Lemmas 2.4 and 3.6, we have 
$\eta\in\Lambda^+$ and $\Omega\Xi(m_1)= g(\eta)m_1$.  Therefore 
$g(\eta)=(rs^{-1})^c$.  \m We now have $g(\mu)=g(\eta)$, where 
$\eta,\mu\in\Lambda^+$, and $\eta\geq \mu$ by construction.  By Lemma 
3.5, $\eta=\mu$, so $M_1$ is the one-dimensional space spanned by $m$, 
and $e_i.m=0 \ (1\leq i<n)$, that is $m\in P$.  This implies $m''=0$, 
a contradiction to the assumption that $M''\neq 0$.  Therefore 
$M''=0$, and $M=M'$ is completely reducible.  \m Finally, we consider 
the case that $M$ does not have weights in $\Lambda$.  We may assume 
that $M$ is indecomposable.  By Remark 2.17, $M$ has all its weights 
in $\chi\cdot\hat{\Lambda}$ for some $\chi$ satisfying (2.15).  By 
Lemma 2.18, $M\cong L_{\chi}\otimes N$ for some $\widetilde{U}$-module 
$N$ whose weights are all in $\Lambda$.  Note that $\widetilde{U}^0$ 
acts semisimply on $N$ as well ($N=L_{\chi^{-1}}\otimes M$), and so 
$N$ is completely reducible by the above argument.  As the tensor 
product of modules distributes over direct sums, $M$ is itself 
completely reducible.  \qed\enddemo \b \n {\bf Remark 3.9.} It is 
necessary to include the hypothesis that $\widetilde{U}^0$ acts 
semisimply in Theorem 3.8, as the next examples illustrate.  (Recall 
that $\widetilde{U}^0$ does indeed act semisimply on any simple 
$\widetilde{U}$-module, as remarked in the text following (1.6).) Let 
$V=\K^m$ for $m\geq 2$ and $\xi,\xi'\in \K\setminus\{0\}$.  We define 
a $\widetilde{U}$-module structure on $V$ by requiring that $e_i,f_i$ 
act as multiplication by 0 and $a_i, b_i$ act via the $m\times m$ 
Jordan blocks with diagonal entries $\xi,\xi'$, respectively.  The 
relations of $\widetilde{U}$ hold on $V$: (R1) is satisfied as these 
matrices are invertible and commute with one another.  (R4) holds as 
$a_ib_{i+1}$ and $a_{i+1}b_i$ act via the same matrix.  The remaining 
relations hold as $e_i,f_i$ act as multiplication by 0.  The scalars 
$\xi,\xi'$ may be chosen so that $V$ has weights in $\Lambda$, for 
example choose an integer $c$, let $\lambda=c(\epsilon_1 +\cdots 
\epsilon_n)$, and set $\xi=r^c=\hat{\l}(a_i), \ 
\xi'=s^c=\hat{\l}(b_i)$.  Clearly $V$ is not completely reducible as 
the Jordan blocks are not diagonalizable.  A related example for 
$U_{r,s}({\frak {sl}}_n)$ is given by sending $\omega_i, \omega_i'$ to 
the same Jordan block with diagonal entries $\xi_i\in\K-\{0\} \ (1\leq 
i< n)$.

\b \m \head \S 4.  The $R$-matrix 
\endhead \m
 In this section we recall the definition of the $R$-matrix 
from [BW, Sec.  4] and use it to prove a more general result on 
commutativity of the tensor product of finite-dimensional
modules than was given there
(compare [BW, Thm. 4.11] with Theorem 4.2 below).  Let $M,M'$ be 
$\widetilde{U}$-modules in category $\Cal O$.  We define an 
isomorphism of $\widetilde{U}$-modules $R_{M',M}: M'\otimes M 
\rightarrow M\otimes M'$ as follows.
If $\lambda=\sum_{i=1}^n\lambda_i\alpha_i\in\Lambda$, where 
$\alpha_n=\epsilon_n$, set 
$$\aligned
  \omega_{\l}&=\omega_1^{\l_1}\cdots \omega_{n-1}^{\l_{n-1}}a_n^{\l_n}\\
  \omega'_{\l}&=(\omega'_1)^{\l_1}\cdots(\omega'_{n-1})^{\l_{n-1}}
   b_n^{\l_n}. \endaligned $$
Also let 
$$\Theta =\sum_{\zeta\in Q^+}\sum_{k=1}^{d_{\zeta}} v^{\zeta}_k\otimes
u^{\zeta}_k,$$
where the notation is as in the paragraph following (3.1). Define
$$R_{M',M}=\Theta\circ \widetilde{f}\circ P$$
where $P(m'\otimes m)=m\otimes m'$, $\widetilde{f}(m\otimes m')=
(\omega_{\mu}',\omega_{\l})^{-1}(m\otimes m')$ when $m\in M_{\l}$
and $m'\in M'_{\mu}$, and the Hopf pairing $( \ , \ )$ 
is defined in (3.1). 
Then $R_{M',M}$ is an isomorphism of $\widetilde{U}$-modules that
satisfies the quantum Yang-Baxter equation and the hexagon identities
[BW, Thms. 4.11, 5.4, and 5.7].
\m

We will show that the tensor product of {\it any} 
two finite-dimensional $\widetilde{U}$-modules in $\Cal O$
is commutative (up to module isomorphism).
We first prove this in the special case that one of the modules is a 
one-dimensional module $L_{\chi}=L(\chi)$, as defined in Section 2. 
 \b 
\proclaim{Lemma 4.1} Let $M$ be a
$\widetilde{U}$-module in category $\Cal O$, and let $L_{\chi}$ 
be a one-dimensional $\widetilde{U}$-module.  Then $$L_\chi\otimes M 
\cong M\otimes L_\chi.$$\endproclaim \m \demo{Proof} Fix a basis 
element $v$ of $L_\chi$.  Define a linear function $F:L_\chi\otimes 
M\rightarrow M\otimes L_\chi$ as follows.  If $m\in M_{\lambda}$, 
where $\lambda =-\sum_{i=1}^{n} c_i\alpha_i$, then $$F(v\otimes m)= 
\chi_1^{c_1}\cdots \chi^{c_n}_n
m\otimes v,$$
where $\chi_i=\chi(\w_i)=\chi(\w_i') \ (1\leq i<n)$ and $\chi_n =\chi(a_n)$.
Clearly $F$ is bijective, and we check that $F$ is a 
$\widetilde{U}$-homomorphism:
$$\aligned
e_i.F(v\otimes m) &= \chi_1^{c_1}\cdots \chi^{c_n}_n 
(e_i\otimes 1 +\omega_i\otimes e_i)(m\otimes v)\\
 &= \chi_1^{c_1}\cdots\chi^{c_n}_n
e_i.m\otimes v.\endaligned$$
On the other hand, as $e_i.m\in M_{\lambda +\alpha_i}$, we have
$$\aligned
F(e_i.(v\otimes m)) &= F((e_i\otimes 1 +\omega_i\otimes e_i)(v\otimes m))\\
  &= \chi_i F(v\otimes e_i m)\\
  &= \chi_i (\chi_1^{c_1}\cdots \chi_i^{c_i-1}\cdots
            \chi^{c_n}_n) e_i.m\otimes v\\
  &= e_i.F(v\otimes m).\endaligned$$
Similarly, $F$ commutes with $f_i$.
As the action by $a_i,b_i$ preserves the weight spaces, $F$ commutes with
$a_i, b_i \ (1\leq i\leq n)$ as well.  Therefore $F$ is an isomorphism of
$\widetilde{U}$-modules. \qed\enddemo
\b
\proclaim{Theorem 4.2} Let $M$ and $M'$ be finite-dimensional
$\widetilde{U}$-modules with $\widetilde U^0$ acting semisimply.  Then
$$M\otimes M'\cong M'\otimes M.$$
\endproclaim
\m
\demo{Proof} As the tensor product distributes over direct sums, we may assume
that $M$ and $M'$ are indecomposable.  Therefore the weights of $M$
are all in $\chi\cdot\hat{\Lambda}$ for some algebra homomorphism 
$\chi:\widetilde{U}^0\rightarrow\K$ with
$\chi(\w_i) = \chi(\w_i')$.  (See Remark 2.17.)   
By Lemma 2.18, $M\cong L_{\chi}\otimes 
N$ for some module $N$ with weights in $\Lambda$.  Similarly, $M'\cong L_{\chi'}\otimes
N'$ for some $\chi'$.  By Lemma 4.1 and [BW, Thm. 4.11], $$\aligned 
M\otimes M' \cong (L_{\chi}\otimes N)\otimes (L_{\chi'}\otimes N') 
&\cong (L_{\chi}\otimes L_{\chi'})\otimes (N \otimes N')\\
  &\cong(L_{\chi'}\otimes L_{\chi})\otimes (N' \otimes N)\\
  &\cong (L_{\chi'}\otimes N')\otimes (L_{\chi}\otimes N)
  \cong M'\otimes M.\qed\endaligned$$ \enddemo
\b
\m
\head \S 5. Tensor powers of the natural module \endhead
\m
In this section we consider tensor powers $V^{\ot k} =
V \ot V \ot \cdots \ot V$ ($k$ factors) of the natural
module $V$ for $\widetilde U$ (defined in Section 1).  
Set $R = R_{V,V}$, and for
$1 \leq i < k$, let $R_i$ be the $\widetilde
U$-module isomorphism on $V^{\ot k}$ defined by

$$R_i(z_1 \ot z_2 \ot \cdots \ot z_k)
= z_1 \ot \cdots \ot z_{i-1}\ot R(z_i \ot z_{i+1}) \ot z_{i+2} \ot \cdots \ot
z_k.$$

\n Then it is a consequence of the quantum Yang-Baxter equation
that the braid relations 

$$\aligned  R_i \circ  R_{i+1} \circ  R_i
& =  R_{i+1} \circ  R_i \circ  R_{i+1} \quad \quad \text{for} \quad 1 \leq i < k \\
 R_i \circ  R_j & =  R_j \circ  R_i 
\qquad \qquad \text{for} \quad |i-j| \geq 2 \endaligned
\tag 5.1$$

\n hold.
\m
We would like to argue that

$$ R_i^2 = (1-rs^{-1}) R_i + rs^{-1}\id \tag 5.2$$

\n for all $i=1,\dots, k-1$.   For this it suffices to work with the
2-fold tensor product $V \ot V$. 
\b  

\proclaim{Proposition 5.3} Whenever $s \neq -r$, the module
$V \ot V$ decomposes into two simple submodules, $S^2_{r,s}(V)$
(the $(r,s)$-symmetric tensors) and $\wedge^2_{r,s}(V)$
(the $(r,s)$-antisymmetric tensors).  These modules
are defined as follows:
\roster
\item"{(i)}" $S^2_{r,s}(V)$ is the span of
$\{ v_i \ot v_i \mid 1\leq i \leq n\} \cup \{v_i \ot v_j +
s v_j \ot v_i \mid 1 \leq i < j \leq n\}$.
\m 
\item"{(ii)}" $\Lambda^2_{r,s}(V)$ is the span of
$\{v_i \ot v_j -rv_j \ot v_i \mid 1 \leq i < j \leq n\}.$
\endroster
\endproclaim 
\m
\demo{Proof} The following computations can be used to see that 
$S^2_{r,s}(V)$ and $\Lambda^2_{r,s}(V)$ are submodules:

$$\aligned
e_k.(v_i \ot v_i) &= \delta_{i,k+1}(v_k \ot v_{k+1} + s v_{k+1} \ot v_k)\\
f_k.(v_i\ot v_i) &= \delta_{i,k}(v_k\ot v_{k+1} + sv_{k+1}\ot v_k)
\endaligned $$
$$\aligned
e_k.(v_i \ot v_j + s v_j \ot v_i)
& = \cases 0 &\quad \  \qquad \text{if} \ \ k+1 \neq i,j \\
v_k \ot v_j + s v_j \ot v_k & \quad \  \qquad \text{if} \ \ k+1 = i \\
v_i \ot v_k + s v_k \ot v_i & \quad \ \qquad \text{if} \ \ k+1 = j, \ k \neq i\\
(r+s) v_k \ot v_k &\quad \ \qquad \text{if} \ \ k+1 = j, \ k = i \endcases 
\endaligned$$
$$\aligned
f_k.(v_i \ot v_j + s v_j \ot v_i) & = \cases 0 & \qquad \text{if} \ \ k 
\neq i,j \\
v_i \ot v_{k+1} + s v_{k+1} \ot v_i & \qquad \text{if} \ \ k = j \\
v_{k+1} \ot v_j + s v_j \ot v_{k+1} & \qquad \text{if} \ \ k = i, \ k+1\neq j\\
(r+s) v_{k+1} \ot v_{k+1} & \qquad \text{if} \ \ k = i, \ k +1= j \endcases 
\endaligned$$
$$\aligned
e_k.(v_i \ot v_j -r v_j \ot v_i)
& = \cases 0 & \quad \ \qquad \text{if} \ \ k+1 \neq i,j \\
v_k \ot v_{j} -r v_{j} \ot v_k & \quad \  \qquad \text{if} \ \ k+1=i\\
v_{i} \ot v_k -r v_k \ot v_{i} & \quad \  \qquad \text{if} \ \ k +1 = j, \ k 
\neq i \\
0 & \quad \ \qquad \text{if} \ \ k+1 = j, \ k = i\endcases
\endaligned$$
$$\aligned
f_k.(v_i \ot v_j -r v_j \ot v_i)
& = \cases 0 & \qquad \text{if} \ \ k \neq i,j \\ 
v_i \ot v_{k+1} -r v_{k+1} \ot v_i & \qquad \text{if} \ \ k=j\\
v_{k+1} \ot v_j -r  v_j \ot v_{k+1} & \qquad \text{if} \ \ k  = i, \ k+1 \neq j \\
0 & \qquad \text{if} \ \ k = i, \ k+1 = j.\endcases
\endaligned$$

Note that each weight space of $S^2_{r,s}(V)$
is one-dimensional and is spanned by one of the weight vectors listed in (i). 
Therefore any submodule of $S^2_{r,s}(V)$ must
contain one of these vectors.  The above computations show that any of
these vectors generates all of $S^2_{r,s}(V)$ in case $s\neq -r$.  In 
particular, $v_1\ot v_1$ is a highest weight vector, and it is easy to see
that given any other vector in (i), there is an element of $U$ taking it
to $v_1\ot v_1$. Therefore $S^2_{r,s}(V)$ is simple.  A similar
argument proves that $\Lambda^2_{r,s}(V)$ is simple, with highest weight 
vector $v_1\ot v_2-rv_2\ot v_1$.  \qed \enddemo \b \n {\bf Remark 
5.4.} The $s = -r$ case is ``nongeneric,'' and in this exceptional 
case, $V\ot V$ need not be completely reducible.  For example, when $n 
= 2$ what happens is that $v_1 \ot v_2 - r v_2 \ot v_1$ spans a 
one-dimensional module (as it does for $n = 2$ generic) that is not 
complemented in $V\ot V$.  Modulo that submodule, $v_1 \ot v_1$ spans 
a one-dimensional module.  Modulo the resulting two-dimensional 
module, $v_1 \ot v_2 + r v_2 \ot v_1$ and $v_2 \ot v_2$ span a 
two-dimensional module.

\b \proclaim {Proposition 5.5} The 
minimum polynomial of $R = R_{V,V}$ on $V 
\ot V$ is $(t-1)(t+rs^{-1})$ if $s \neq -r$.  \endproclaim 
\m \demo {Proof} It follows from the definition of $ R$ that $ R(v_1 
\ot v_1) = v_1 \ot v_1$ and $ R(v_1 \ot v_2 - r v_2 \ot v_1) = 
-rs^{-1}(v_1 \ot v_2 -rv_2 \ot v_1)$.  By Proposition 5.3, 
$S^2_{r,s}(V)$ and $\Lambda^2_{r,s}(V)$ are simple, and in fact, $v_1 
\ot v_1$ and $v_1 \ot v_2 - r v_2 \ot v_1$ are the highest weight 
vectors.  In particular, each is a cyclic module generated by its 
highest weight vector.  As $ R a(v_1 \ot v_1) = a R(v_1 \ot v_1) = 
a(v_1 \ot v_1)$ for all $a \in \widetilde U$, this implies that 
$S^2_{r,s}(V)$ is in the eigenspace of $R$ corresponding to eigenvalue 
1.  Analogously, $\wedge^2_{r,s}(V)$ corresponds to the eigenvalue 
$-rs^{-1}$, and since $V \ot V$ is the direct sum of those submodules, 
we have the desired result.  \qed \enddemo
 
 \b From Proposition 5.5 it follows that $R$ acts as 
 $$r\sum_{i<j}E_{j,i}\otimes E_{i,j} + s^{-1}\sum_{i<j}E_{i,j}\otimes 
 E_{j,i} + (1-rs^{-1})\sum_{i<j}E_{j,j}\otimes E_{i,i} 
 +\sum_iE_{i,i}\otimes E_{i,i}
\tag 5.6 $$
on $V\otimes V$.  Indeed, (5.6) is a linear operator that acts on
$S^2_{r,s}(V)$ as multiplication by 1, and on $\Lambda^2_{r,s}(V)$ as
multiplication by $-rs^{-1}$.  By Proposition 5.5, $R$ has the same
properties, and so $R$ is equal to this sum on $V\otimes V$.
\b
\m
\head \S 6.  Quantum Schur-Weyl duality \endhead
\m 

Assume $r,s \in \K$.
Let $H_k(r,s)$ be the unital associative algebra over $\K$ with
generators $T_i$, \ $1 \leq i < k$, subject to the relations:
\m
\roster
\item"{(H1)}" $T_i T_{i+1}T_i  =  T_{i+1}T_iT_{i+1}, \ \ \qquad 1 \leq i
< k$  
\item"{(H2)}" $T_i T_j = T_j T_i$, \ \ \qquad $|i-j| \geq 2$
\item"{(H3)}" $T_i^2 = (s-r)T_i + rs 1$.
\endroster
\m
When $r \neq 0$, the elements  $t_i = r^{-1}T_i$ satisfy the braid relations (H1),
(H2), along with the relation 

$$t_i^2 = (q-1)t_i + q 1,  \tag H3'$$

\n  where $q = r^{-1}s$.  The {\it Hecke algebra} $H_k(q)$ (of
type A$_{k-1}$)  
is generated by elements $t_i$, $1 \leq i < k$, which satisfy 
(H1), (H2), (H3'). 
It has dimension $k!$ and is semisimple whenever
$q$ is not a root of unity.   At $q = 1$, the Hecke algebra $H_k(q)$ is isomorphic to
$\K
S_k$, the group algebra of the symmetric group $S_k$, where we may identify $t_i$ with
the transposition $(i \  i+1)$.   
\m
The two-parameter Hecke algebra $H_k(r,s)$ defined above is isomorphic
to $H_k(r^{-1}s)$
whenever
$r
\neq 0$.  Thus, it is semisimple whenever $r^{-1}s$ is not a root of unity.  For any $\sigma
\in S_k$, we may define
$T_{\sigma} = T_{i_1}
\cdots T_{i_\ell}$ where $\sigma = t_{i_1}\dots t_{i_\ell}$ is a reduced expression for $\sigma$
as a product of transpositions.  It follows from (H1) and (H2) that $T_{\sigma}$ is independent
of the reduced expression and these elements give a basis.  
\m
The results of Section 5 show that the $\widetilde U$-module $V^{\ot k}$
affords a representation of the Hecke algebra $H_k(r,s)$:

$$\aligned H_k(r,s) & \rightarrow \End_{\widetilde U}(V^{\ot k}) \\
T_i & \mapsto s R_i \qquad (1 \leq i < k). \endaligned \tag 6.1$$

When $k = 2$ and  $s \neq -r$,  
$V^{\ot 2} =  S^2_{r,s}(V) \oplus \wedge^2_{r,s}(V)$ is a decomposition
of $V^{\ot 2}$ into simple $\widetilde U$-modules by Proposition 5.3.  The 
maps $p_1 = (sR_1+r)/(s+r)$ and $p_2 = (s-sR_1)/(s+r)$, $(R_1 = 
R_{V,V})$, are the corresponding projections onto the simple summands.  
Thus, the map in (6.1) is an isomorphism for $k = 2$.  More generally, 
we will show next that it is surjective whenever $rs^{-1}$ is not a 
root of unity, and it is an isomorphism when $n \geq k$.  This is the 
two-parameter version of the well-known result of Jimbo [Ji] that 
$H_k(q) \cong \End_{U_q(\fgl)}(V^{\ot k})$ and is the analogue of 
classical Schur-Weyl duality, $\K S_k \cong \End_{\fgl}(V^{\ot k})$ 
for $n \geq k$.  It requires the following lemma.  The case $n<k$ is 
dealt with separately, and uses the isomorphism $H_k(r,s)\cong 
\End_{\widetilde{U}}(V^{\otimes k})$ of the case $n=k$.  \b 
\proclaim{Lemma 6.2} If $n\geq k$ and $V$ is the natural 
representation of $\widetilde{U}$, then $V^{\otimes k}$ is a cyclic 
$\widetilde{U}$-module generated by $v_1\otimes \cdots\otimes v_k$.  
\endproclaim \m \demo{Proof} Let $\underline{v}=v_1\otimes 
\cdots\otimes v_k$.  We begin by showing that if $\sigma\in S_k$, then 
$v_{\sigma(1)}\otimes\cdots\otimes v_{\sigma(k)}\in 
\widetilde{U}.\underline{v}$.  \m Suppose we have an arbitrary 
permutation $x_1 \ot \cdots \ot x_k$ ($x_i \in \{v_1,\dots,v_k\}$ for 
all $i$) of the factors of $\underline{v}$.  For some $\ell < m$, 
assume that $x_\ell= v_j$ and $x_m = v_{j+1}$.  Then because of the 
formulas

$$\aligned 
\Delta^{k-1}(e_j)&= \sum_{i=1}^{k}
\underbrace{\omega_j\otimes\cdots\ot\omega_j}_{i-1}\otimes e_j \otimes 
\underbrace{1\otimes\cdots\otimes 1}_{k-i}€ \\ 
\Delta^{k-1}(f_j)&=\sum_{i=1}^{k}
\underbrace{1\otimes\cdots\otimes 1}_{k-i}\otimes f_j \otimes
\underbrace{\omega_j'\otimes\cdots\otimes\omega_j'}_{i-1}, \endaligned \tag 6.3 
$$ there are nonzero scalars $c$ and $c'$ such that 
$$(ce_jf_j+c').(x_1\otimes\cdots\otimes x_k)
=x_1\otimes\cdots\otimes x_m
\otimes\cdots\otimes x_\ell \otimes \cdots\otimes x_k.$$
Similarly, there are nonzero scalars $d$ and $d'$ such that
$$(de_jf_j +d').(x_1\otimes\cdots\otimes x_m
\otimes\cdots\otimes x_\ell \otimes \cdots\otimes 
x_k)=x_1\otimes\cdots\otimes x_k.$$
As the transpositions $(j\ j+1)$ generate $S_k$, $v_{\sigma(1)}\otimes\cdots
\otimes v_{\sigma(k)}\in\widetilde{U}.\underline{v}$ for all $\sigma\in S_k$. 
\m

Next we will use induction on $k$ to establish the following.  For any $k$
elements $i_1,\ldots,i_k\in \{1,\ldots,n\}$ satisfying $i_1\leq i_2\leq
\cdots\leq i_k,$ there is a $u\in \widetilde{U}$ such that
$u.\underline{v}=v_{i_1}\otimes\cdots\otimes v_{i_k}$ and $u$ does not
contain any terms with factors of $e_m,e_{m+1},\ldots, e_{n-1}, f_{m+1},
f_{m+2},\ldots,f_{n-2},$ or $f_{n-1}$ where $m=\text{max}\{i_k,k\}$.
If $k=1$, we may apply $f_{m-1}\cdots f_1$ to $\underline{v}=v_1$ to obtain
$v_m$ for any $m\in \{1,\ldots,n\}$.  
If $k>1$, let $\ell$ be such that $i_{\ell}<i_k$, $i_{\ell +1}=i_{\ell +2}
=\cdots =i_k$.  (If no such $\ell$ exists, that is if $i_1=\cdots =i_k$,
then set $\ell=0$ and apply $u'$ from (6.5) below to $v_1\otimes\cdots \otimes
v_k$ to obtain a nonzero scalar multiple of $v_{i_1}\otimes\cdots \otimes
v_{i_k}$.)  By induction, there is an element $u\in \widetilde{U}$ such
that 
$$u.(v_1\otimes \cdots\otimes v_{\ell})=v_{i_1}\otimes\cdots\otimes
v_{i_{\ell}}, \tag 6.4$$
where $u$ has no terms with factors of $e_{m'},e_{m'+1},\ldots, e_{n-1},
f_{m'+1},\ldots,f_{n-1}$ ($m'=\text{max}\{i_{\ell},\ell\}$).
\m
Suppose initially that $i_{\ell}\leq \ell$.  Then $m'=\ell$, and so 
$u.(v_1\otimes\cdots\otimes v_k)$ is a nonzero scalar multiple of
$(v_{i_1}\otimes \cdots\otimes v_{i_{\ell}})\otimes (v_{\ell +1}\otimes\cdots
\otimes v_k)$. Now apply

$$u' =  
\cases(f_{i_k -1}f_{i_k -2}\cdot \cdot f_{\ell +1})\cdots (f_{i_k -1}f_{i_k 
-2})(f_{i_k -1}) (e_{i_k}e_{i_k+1}\cdot \cdot e_{k-1})\cdots 
(e_{i_k}e_{i_k+1})(e_{i_k}) \\
 \hskip 4.2 truein  \text{ if } \ i_k < k  \\
(f_{i_k-1}f_{i_k-2}\cdots f_{\ell +1})\cdots (f_{i_k-1}f_{i_k-2}\cdots
 f_{k-1})(f_{i_k-1}f_{i_k-2}\cdots f_k) \ \ \text{ if } \ i_k\geq k\endcases
 \tag 6.5$$
 
\n to obtain a nonzero scalar multiple of $v_{i_1}\otimes \cdots\otimes 
v_{i_k}$, as desired.  (Note that we did not use any factors of 
$e_m,e_{m+1},\ldots, e_{n-1},f_{m+1},\ldots, f_{n-1}$ 
for $m=\text{max}\{i_k,k\}$.) \m If on the other hand, $i_{\ell}>\ell$ (so 
that $m'=i_{\ell}$ and $i_k>\ell +1$), first apply $u'$ from (6.5) to 
$v_1\otimes \cdots\otimes v_k$ to obtain a nonzero scalar multiple of 
$$(v_1\otimes\cdots\otimes v_{\ell})\otimes 
(v_{i_k}\otimes\cdots\otimes
v_{i_k}),$$
and then apply $u$ from (6.4) to obtain a nonzero scalar multiple of 
$v_{i_1}\otimes\cdots\otimes v_{i_k}$, as desired.
\m Finally, if $i_1,\ldots,i_k\in\{1,\ldots,n\}$ are {\it any} $k$ 
elements (not necessarily in nondecreasing numerical order), let 
$\sigma\in S_k$ be a permutation such that
$$i_{\sigma(1)}\leq i_{\sigma(2)}\leq \cdots\leq i_{\sigma(k)}.$$
By the first paragraph of the proof, there is an element of $\widetilde{U}$
taking $\underline{v}$ to $v_{\sigma^{-1}(1)}\otimes\cdots\otimes 
v_{\sigma^{-1}(k)}$.  Now we may apply $u$ from (6.4) and $u'$ from (6.5)
in the appropriate order (as above) to $v_{\sigma^{-1}(1)}
\otimes\cdots\otimes v_{\sigma^{-1}(k)}$ to obtain a nonzero scalar multiple of
$v_{i_1}\otimes\cdots\otimes v_{i_k}$.
\qed\enddemo
\b
This leads to the two-parameter analogue of Schur-Weyl duality.  \b 
\proclaim{Theorem 6.6} Assume $rs^{-1}$ is not a root of unity.  Then: 
\roster \item"{(i)}" $H_k(r,s)$ maps surjectively onto 
$\End_{\widetilde{U}}(V^{\otimes k})$.

\item"{(ii)}" If $n\geq k$, then $H_k(r,s)$ is isomorphic to 
$\End_{\widetilde{U}}(V^{\otimes k})$.
\endroster
 \endproclaim \m

\demo{Proof} We establish part (ii) first.  Assume $F\in 
\End_{\widetilde{U}}(V^{\otimes k})$ and $\underline{v}=v_1\otimes 
\cdots\otimes v_k$.  As $F$ commutes with the action of 
$\widetilde{U}$, $F(\underline{v})$ must have the same weight as 
$\underline{v}$, that is, $\epsilon_1 +\cdots +\epsilon_k$.  The only 
vectors of $V^{\otimes k}$ with this weight are the linear 
combinations of the permutations of $v_1\otimes\cdots\otimes v_k$, so 
that $$F(\underline{v})=\sum_{\sigma\in S_k}
c_{\sigma}v_{\sigma(1)}\otimes\cdots \otimes v_{\sigma(k)}, \tag 6.7$$ 
for some scalars $c_{\sigma}\in \K$.  
 We will show that there is an element $R^{\sigma}$ in the image of 
$H_k(r,s)$ in $\End_{\widetilde{U}}(V^{\otimes k})$ such that 
$R^{\sigma}(\underline{v})=  
v_{\sigma(1)}\otimes\cdots\otimes v_{\sigma(k)}$.  (Previously
we constructed an element $u \in \widetilde U$ with this property.)
 \m We 
begin with the transposition $\tau = t_j=(j \ j+1)$.  For any tensor 
product $v_{i_1}\otimes \cdots\otimes v_{i_k}$ of distinct basis vectors, we 
have by (5.6) that
 
$$  
v_{i_{\tau (1)}}\otimes \cdots\otimes v_{i_{\tau (k)}}  = \cases r^{-1}R_j 
(v_{i_1}\otimes\cdots \otimes v_{i_k}) & \qquad \text{if \ \ $i_j 
<i_{j+1}$} \\ (sR_j + (r-s)\id) (v_{i_1}\otimes\cdots\otimes v_{i_k}) 
& \qquad \text {if \ \ $i_j >
i_{j+1}$.} \endcases$$

Therefore, if
$\sigma = t_{j_1}\cdots t_{j_m}$, a product of such transpositions, we can
set  $R^{t_{j_{\ell}}}:=r^{-1} R_{j_{\ell}}$ or 
$R^{t_{j_{\ell}}}:=sR_{j_{\ell}} +(r-s)\id$, depending on the numerical 
order of the appropriate indices $i_{j_{\ell}}$ and $i_{j_{\ell}+1}$
in $R^{t_{j_{\ell-1}}}\circ\cdots\circ R^{t_{j_1}} \underline {v}$.  Then 
defining $R^{\sigma}=R^{t_{j_m}}\circ\cdots\circ R^{t_{j_1}} \in 
\End_{\widetilde{U}}(V^{\otimes k})$, we have the desired map such 
that $R^{\sigma}(\underline{v}) = 
v_{\sigma(1)}\otimes\cdots\otimes v_{\sigma(k)}$.  \m Now let 
$F_0=F-\sum_{\sigma\in S_k} c_{\sigma} R^{\sigma}\in \End_ 
{\widetilde{U}}(V^{\otimes k})$ (with the $c_{\sigma}$ coming from 
(6.7)), and note that $F_0(\underline{v})=0$.  As $F_0$ commutes with 
the action of $\widetilde{U}$, we have 
$F_0(\widetilde{U}.\underline{v})=\widetilde{U}.F_0(\underline{v})=0$.  
By Lemma 6.2, $\widetilde{U}.\underline{v}=V^{\otimes k}$.  Therefore 
$F_0$ is the 0-map, which implies $F=\sum_{\sigma\in S_k} c_{\sigma} 
R^{\sigma}$ is in the image of $H_k(r,s)$.  Consequently, the map 
$H_k(r,s) \rightarrow \End_ {\widetilde{U}}(V^{\otimes k})$ in (6.1) 
is a surjection, and  $\End_{\widetilde{U}}(V^{\otimes k})$ is the 
$\K$-linear span of $\{R^{\sigma}\mid \sigma\in S_k\}$.  Now suppose 
that $\sum_{\sigma\in S_k}c_{\sigma}R^{\sigma} =0$ for some scalars 
$c_{\sigma}\in \K$.  Then in particular, $$0=\sum_{\sigma\in 
S_k}c_{\sigma}R^{\sigma}(\underline{v})=\sum_
{\sigma\in S_k}c_{\sigma} v_{\sigma(1)}\otimes \cdots\otimes v_{\sigma(k)}.$$
The vectors $\{v_{\sigma(1)}\otimes\cdots\otimes v_{\sigma(k)}\mid\sigma\in
S_k\}$ are linearly independent, so $c_{\sigma}=0$ for all $\sigma\in S_k$.
This implies that $\{R^{\sigma}\mid \sigma\in S_k\}$ is a basis for the vector
space $\End_{\widetilde{U}}(V^{\otimes k})$ and
$\dim_{\K}\End_{\widetilde{U}}(V^{\otimes k}) = k!$\,.  Because $H_k(r,s)$ is 
isomorphic to $H_k(r^{-1}s)$, it has dimension $k!$ also.  Therefore, 
$H_k(r,s)$ is isomorphic to $\End_{\widetilde{U}}(V^{\otimes k})$ for 
$n \geq k$, as asserted.

 \m Next we turn to the 
proof of (i) and assume here that $n < k$.  For $i=n,k$, let 
$\widetilde{U}_i=U_{r,s}({\frak {gl}}_i)$, let $\Lambda_i$ be the 
weight lattice of ${\frak {gl}}_i$, and let $V_i$ be the natural 
$\widetilde{U}_i$-module.  By (ii), we may identify $H_k(r,s)$
with $\End_{\widetilde{U}_k}(V_k^{\otimes k})$.  We will show that 
$H_k(r,s)$ maps surjectively onto End$_{\widetilde{U}_n}(V_n^{\otimes 
k})$.  \m Consider $V_k^{\otimes k}$ as a $\widetilde{U}_n$-module via 
the inclusion of $\widetilde{U}_n$ into $\widetilde{U}_k$, and regard 
$V_n^{\otimes k}$ as a $\widetilde{U}_n$-submodule of $V_k^{\otimes 
k}$ in the obvious way.  Now $V_n^{\otimes k}$ is a finite-dimensional 
$\widetilde{U}_n$-module on which $\widetilde{U}_n^0$ acts semisimply, 
so by Theorem 3.8, it is completely reducible.  Therefore,
$$V_n^{\otimes k} =  L_1\oplus\cdots\oplus L_t \tag 6.8 $$
for simple $\widetilde{U}_n$-modules $L_i$.  It suffices
to show that the projections onto the simple summands $L_i$
can be obtained from $H_k(r,s)$.  
\m Consider $$\widetilde{U}_k .  V_n^{\otimes k} = \widetilde{U}_k 
.  L_1
+\cdots + \widetilde{U}_k . L_t, \tag 6.9 $$
the $\widetilde{U}_k$-submodule of $V_k^{\otimes k}$ generated by 
$V_n^{\otimes k}$.  By Corollary 2.5, each $L_i$ is isomorphic to 
some $L(\lambda _i)$, $\lambda_i\in \Lambda_n ^+$, and in particular 
is generated by a highest weight vector $m_i$ with  $e_j .  m_i=0$ for all 
$j$ ($1\leq j <n$).  We claim that $e_j .  m_i=0$ as well when $n\leq 
j<k$.  This follows from the expression for $\Delta^{k-1}(e_j)$ in 
(6.3) and the action of $e_j$ on the natural module $V_k$ for 
$\widetilde{U}_k$ given by $e_j.v_i=\delta_{i,j+1}v_j$, because $m_i$ 
must be some linear combination of vectors $v_{i_1}\otimes 
\cdots\otimes v_{i_k}$ with $i_1,\ldots,i_k\in\{1,\ldots, n\}$.  
Therefore $m_i$ is also a highest weight vector for the 
finite-dimensional $\widetilde{U}_k$-module $\widetilde{U}_k .  L_i$.  
By Theorem 2.1 and Lemma 3.7, $\widetilde{U}_k .  L_i= \widetilde{U}_k 
.  m_i$ is a simple $\widetilde{U}_k$-module.  Therefore (6.9) must be 
a direct sum:

 $$\widetilde{U}_k .  V_n^{\otimes k} = \widetilde{U}_k 
.  L_1\oplus
\cdots\oplus\widetilde{U}_k . L_t.$$

Because $V_k^{\otimes k}$ is a completely reducible 
$\widetilde{U}_k$-module, there is some complementary 
$\widetilde{U}_k$-submodule $W$ such that $$V_k^{\otimes k} = \widetilde{U}_k 
.  L_1 \oplus\cdots\oplus\widetilde{U}
_k . L_t \oplus W. \tag 6.10$$
Let $\pi_i \in H_k(r,s)$ be the projection of $V_k^{\otimes k}$ onto 
$\widetilde{U}_k .  L_i$.  Then, $\pi_i$ commutes
with the $\widetilde U_k$-action, and  acts as the identity map on 
$\widetilde{U}_k .  L_i$ and as $0$ on the other summands in (6.10).  
Since $L_j \subseteq \widetilde{U}_k .  L_j$ for all $j$, the map 
$\pi_i$ restricted to $ V_n^{\otimes k}$ commutes with the $\widetilde 
U_n$-action and is the projection onto $L_i$.  Thus, $H_k(r,s) 
\rightarrow \End_{\widetilde U_n}(V_n^{\otimes k})$ is onto.  \qed 
\enddemo \m \b

\vskip 3 mm
 
\address 
\newline 
Department of Mathematics, University of Wisconsin, Madison,
Wisconsin 53706   \newline 
benkart\@math.wisc.edu
\newline 
 \newline  
Department of Mathematics and Statistics, University of Massachusetts, Amherst,
Massachusetts  01003  
\newline 
(2001-02) Department of Mathematics and Computer Science,
Amherst College, Amherst, Massachusetts 01002
\newline
wither\@math.umass.edu 
\endaddress 

\enddocument
\end


\vskip 3 mm
 
\address 
\newline 
Department of Mathematics, University of Wisconsin, Madison,
Wisconsin 53706   \newline 
benkart\@math.wisc.edu
\newline 
 \newline  
Department of Mathematics, University of Massachusetts, Amherst,
Massachusetts  01003  
\newline
wither\@math.umass.edu 
\endaddress 


\Refs 
\widestnumber\key{KMP} 
\m
\ref \key BR \by G. Benkart and T. Roby \paper Down-up algebras
\jour J. Algebra \vol 209 \yr 1998 \pages 305-344 \, \moreref
\jour Addendum \vol 213 \yr 1999 \pages 378 \endref
\m
\ref \key BW \by G. Benkart and S. Witherspoon \paper Two-parameter quantum
groups and Drinfel'd doubles \jour preprint \endref
\m 
 \ref \key H \by T.  Halverson \book Characters of the Centralizer 
Algebras for Mixed Tensor Representations of the General Linear Group 
and its $q$-Deformations \publ Ph.D.  Thesis, University of Wisconsin 
- Madison, 1993 \endref \m 
\ref \key Ja \by J. C. Jantzen \book Lectures on Quantum Groups \vol 6 
\publ Graduate Studies in Math., Amer. Math. Soc. \yr 1993 \endref
\m
\ref \key Ji \by M.  Jimbo \paper A 
$q$-analogue of $U(\fgl(N+1))$, Hecke algebra, and the Yang-Baxter 
equation, Lett.  Math.  Phys.  \vol 11 \yr 1986 \pages 247-252 \endref

\m \ref \key Jo \by A.  Joseph \book Quantum Groups 
and Their Primitive Ideals \publ Ergebnisse der Mathematik und ihrer 
Grenzgebiete, Springer-Verlag, Berlin \yr 1995 \endref \m \ref \key L 
\by G.  Lusztig \book Introduction to Quantum Groups \publ 
Birkh\"{a}user \yr 1993 \endref \m \ref\key T \by M.  Takeuchi \paper 
A two-parameter quantization of $GL(n)$ \jour Proc.  Japan Acad.  \vol 
66 Ser.  A \yr 1990 \pages 112-114 \endref

\endRefs

\Refs 
\widestnumber\key{KMP} 
\ref \key BL \by V.V. Bavula and T.H. Lenagan  \paper Generalized
Weyl algebras are tensor Krull minimal \jour preprint \endref
\m
\ref\key B \by G. Benkart \paper
Down-up algebras and Witten's deformations of
the universal enveloping algebra of $\fsl$\; \moreref 
\book Recent Progress in Algebra \pages 29-45  \bookinfo
Contemp. Math. \vol 224  \publaddr  Amer. Math. Soc. \yr 1998  \endref
\m
\ref \key BR \by G. Benkart and T. Roby \paper Down-up algebras
\jour J. Algebra, \vol 209 \yr1998 \pages 305-344  \; \moreref
\jour Addendum \vol 213 
\yr 1999 \pages 378 \endref
\m
\ref \key BW \by G. Benkart and S. Witherspoon \paper A Hopf structure
for down-up algebras, Math. Zeitschrift, \toappear  \endref
\m
\ref \key CM \by P.A.A.B. Carvalho and I.M. Musson \paper Down-up
algebras and their representation theory \jour J. Algebra 
\vol 228 \pages 286-310 \yr 2000 \endref 
\m
\ref \key H \by T. Halverson \book Characters of the Centralizer Algebras for
Mixed Tensor Representations of the General Linear Group and its $q$-Deformation
 \publ Ph.D. Thesis, University of Wisconsin - Madison 1993  \endref
\m
\ref \key Ja \by J.C. Jantzen \book Lectures on Quantum Groups \publ Graduate
Studies in Math.,  Amer. Math. Soc. \vol 6  
\yr 1996 \publaddr Providence \endref 
\m
\ref \key Ji \by M. Jimbo \paper A $q$-analogue of
$U(\fgl(N+1))$, Hecke algebra, and the Yang-Baxter equation, Lett. Math.
Phys. \vol 11 \yr 1986 \pages 247-252 \endref
\m 
\ref \key Jin \by N.H. Jing \paper Quantum groups with two parameters
\; \moreref 
\book  Deformation Theory and Quantum
Groups with Applications to Mathematical Physics (Amherst, MA, 1990)  \pages 129-138 
 \bookinfo
Contemp. Math. \vol 134  \publaddr  Amer. Math. Soc. \yr 1992  \endref
\m
\ref 
\key Jo  \by A. Joseph \book
 Quantum Groups and Their Primitive Ideals  \publ Ergebnisse der Mathematik und
ihrer Grenzgebiete, Springer-Verlag \yr 1995 \publaddr Berlin  
\endref   
 \m
\ref \key Jor \by D.A. Jordan \paper
Down-up algebras and ambiskew polynomial rings \jour J. Algebra \vol 228
\pages 311-346 \yr 2000 \endref
\m
\ref \key KK1 \by E. Kirkman and J. Kuzmanovich \paper Primitivity of Noetherian down-up
algebras  \jour Comm. Algebra \vol 28 \pages 2983-2997 \yr 2000\endref
\m
\ref \key KK2 \by E. Kirkman and J. Kuzmanovich \paper Non-Noetherian down-up algebras 
\jour preprint \endref  
\m
\ref \key KMP \by E.E. Kirkman, I. Musson, and D. Passman \paper
Noetherian down-up algebras \jour Proc. Amer. Math. Soc. \vol 127
\yr 1999 \pages 3161-3167 \endref  
\m
\ref\key K \by P.P. Kulish \paper A two-parameter quantum group and gauge
transformations (in Russian) \jour Zap. Nauch. semin. LOMI \vol 180 \yr 1990
\pages 89-93 \endref
\m
\ref \key Ku \by R.S. Kulkarni \paper Down-up algebras
and their representations \jour J. Algebra \toappear \endref  
\m
\ref \key M \by S. Montgomery \book Hopf Algebras and Their Actions on Rings. CBMS
Conf. Math. Publ.  \vol 82  \publ Amer. Math. Soc. 
Providence \yr 1993 \endref
\m 
\ref \key T \by M. Takeuchi  \paper A two-parameter quantization of
GL(n) \jour Proc. Japan Acad. \vol 66  Ser. A \yr 1990
\pages 112-114 \endref
 
\endRefs
\enddocument